\documentclass[journal]{IEEEtran}
\pdfoutput=1
\usepackage{graphicx}
\usepackage{epstopdf}
\usepackage{mathrsfs}
\usepackage{amsfonts}
\usepackage{multirow}

\usepackage{ifpdf}

\usepackage{cite}
\usepackage[numbers,sort&compress]{natbib}

\ifCLASSINFOpdf
\else
\fi
\usepackage[cmex10]{amsmath}
\usepackage[tight,footnotesize]{subfigure}

\usepackage{cases}
\begin{document}
%
\title{Multilevel Fast Multipole Algorithm for Characteristic Mode Analysis}
%
%
%

\author{Qi~I.~Dai,
	 Jun Wei~Wu,
	 Ling Ling~Meng,
        Weng~Cho~Chew,~\IEEEmembership{Fellow,~IEEE,} 
        Wei~E.~I.~Sha
\thanks{This work was supported in part by the National Science
Foundation under Award Number 1218552, in part by the Research Grants
Council of Hong Kong (GRF 711609 and 711508), and in part by the University
Grants Council of Hong Kong (Contract AoE/P-04/08).}
\thanks{Q. I. Dai, L. L. Meng and W. C. Chew are with the Department of Electrical and Computer Engineering, University of Illinois at Urbana$-$Champaign, Urbana, IL 61801 USA.}
\thanks{J. W. Wu is with the School of Electronics and Information, Wuhan University, Hubei, 430072, China.}
\thanks{W. E. I. Sha is with the Department of Electrical and Electronic Engineering, the University of Hong Kong, Hong Kong.}
}

\maketitle

\begin{abstract}
Characteristic mode (CM) analysis poses challenges in computational electromagnetics (CEM) as it calls for efficient solutions of dense generalized eigenvalue problems (GEP). Multilevel fast multipole algorithm (MLFMA) can greatly reduce  the computational complexity and memory cost for matrix-vector product operations, which is powerful in iteratively solving large scattering problems. In this article, we demonstrate that MLFMA can be easily incorporated into the implicit restarted Arnoldi (IRA) method for the calculation of CMs, where MLFMA with the sparse approximate inverse (SAI) preconditioning technique is employed to accelerate the construction of Arnoldi vectors.   This work paves the way of CM analysis for large-scale and complicated three-dimensional ($3$-D) objects with limited computational resources.

\end{abstract}


%
\IEEEpeerreviewmaketitle

\section{Introduction}
%
%
%
%
\IEEEPARstart{A}{fter} its humble beginnings in the 1970s, characteristic mode (CM) theory have gained a recent resurgence of interest in the field of antenna design and optimization. Initiated by Garbacz and refined by Harrington and Mautz \cite{Garbacz, Harrington1971,Harrington1971a}, CM theory was popularized in the antenna community by the work of Cabedo-Fabres as it has been shown promising for systematic antenna design \cite{Fabres2007}. Such a systematic approach becomes increasingly favored since heuristic approaches based on experience and intuition of antenna engineers hardly survive the demanding requirement of designing complex antenna systems. CM analysis is capable of characterizing arbitrary object's radiation and scattering properties only relying on its geometry and material properties rather than source configuration, which makes it a favorable candidate for systematic design methods where complicated functionality tradeoffs can be handled. Moreover, CM analysis provides useful physical insight into antenna operation as well as guidance on the excitation of desired radiation modes \cite{Ethier12, Obeidat10, Adams11}. 
 
Although the use of CM analysis has shown success in antenna design, relatively little effort has been devoted to the numerical computation of CMs. Development of fast CM solvers becomes indispensable for applications such as antenna synthesis where sources are integrated to electrically large platforms, and multiscale modeling that contains fine structures. Current strategy can only deal with small problems as it explicitly generates the impedance matrix by applying a Galerkin's procedure to the electric field integral equation (EFIE). The memory cost is $O(N^2)$ where $N$ is the number of unknowns. When $N$ is in the order of $10^2$ or $10^3$, one can find all eigenpairs (eigenvalues and eigenvectors) of the dense generalized eigenvalue problem (GEP) by Schur decomposition with an overall computational complexity of $O(N^3)$. However, in most instances, only a small portion of the spectrum is desirable, which can be efficiently solved for using Krylov subspace iterations such as Lanczos and Arnoldi algorithms. A number of matrix-vector products are performed in such methods, each of which has a computational complexity of $O(N^2)$. Even so, this $O(N^2)$ complexity hinders the use of CM theory in large object analysis where $N$ may be in the order of $10^5$ or even larger.

The fast multipole method (FMM) was originally proposed by Rokhlin for fast particle simulations where static integral equations are solved \cite{Rokhlin87}. It was later extended for the calculation of  electromagnetic scattering by many research groups \cite{Coifman93,Song94,Song97,Chew2000}. A two-level FMM reduces the complexity of a matrix-vector product and memory cost from $O(N^2)$ to $O(N^{1.5})$. A nonnested two-level FMM incorporated with ray propagation physics further reduces the complexity to $O(N^{4/3})$ \cite{Wagner94}. 
Further reduction in complexity and memory cost was achieved by the development of the multilevel fast multipole algorithm (MLFMA). The complexity becomes $O(N \log^2 N)$ when MLFMA is implemented using signature function, interpolation, and filtering \cite{Dembart95}. Based on translation, interpolation, anterpolation, and a grid-tree data structure, MLFMA reduces the complexity and memory cost to $O(N \log N)$ \cite{Song94,Song97,Chew2000}.  
 
In this paper, MLFMA is successfully adopted in CM analysis. We demonstrated that MLFMA can be easily incorporated into iterative eigensolvers such as the implicit restarted Arnoldi (IRA) method \cite{Sorensen} and Jacobi-Davidson QZ (JDQZ) method \cite{Fokkema}. The sparse approximate inverse (SAI) of the explicitly available submatrix for near group interaction is used as a preconditioner \cite{Jeonghwa} when Arnoldi vectors in IRA are constructed . Other preconditioning schemes such as the Calderon multiplicative preconditioner (CMP) \cite{Andriulli} is briefly discussed, which is promising for CM analysis in even larger-scale applications.

%

\section{Characteristic Mode Formulation}
The widespread CM formulation relies on the EFIE for an arbitrarily shaped perfectly electric conductor (PEC) object given by
\begin{equation}
i k\eta \hat n \times \int_S d\mathbf{r}' \, \overline{\mathbf{G}}(\mathbf{r},\mathbf{r}')\cdot \mathbf {J}(\mathbf{r}') = - 4\pi \hat n\times \mathbf{E}_{inc}(\mathbf{r}), \hspace{.1in} \mathbf{r}\in S
\end{equation}
where $\eta=\sqrt{\mu/\epsilon}$ is the characteristic impedance of free-space, and the Green's dyadic is
\begin{equation}
\overline{\mathbf{G}}(\mathbf{r},\mathbf{r}') = \left \lbrack \overline{\mathbf{I}}-\frac{\nabla \nabla^\prime}{k^2} \right \rbrack 
g(\mathbf{r},\mathbf{r}')
\end{equation}
with
\begin{equation}\label{SGF}
g(\mathbf{r},\mathbf{r}') = \frac{e^{i k  R}}{R } , \qquad R = \vert \mathbf{r} - \mathbf{r}' \vert
\end{equation}
and $\hat n$ is the unit normal of the conducting surface. 

In the method of moments (MoM), the unknown surface current $\mathbf{J}(\mathbf{r})$ is approximated with the popular Rao-Wilton-Glisson (RWG) basis functions $\mathbf{f}_j(\mathbf{r})$ as
\begin{equation}
\mathbf{J}(\mathbf{r}) = \sum_{j=1}^N I_j \mathbf{f}_j(\mathbf{r})
\end{equation}
where $I_j$ are the expansion coefficients. The impedance operator relates the induced current $\mathbf{J}(\mathbf{r})$ to the incident electric field $\hat n \times \mathbf{E}_{inc}(\mathbf{r})$ on $S$, and is projected to form the dense impedance matrix $\bar{\mathbf{Z}}$ using a Galerkin's procedure where testing functions are chosen to be the same as basis functions as
\begin{equation}
\begin{aligned}
Z_{ij} & = i k\eta \left\langle \mathbf{f}_i, \overline{\mathbf{G}}, \mathbf{f}_j \right\rangle_S \\
& = i k\eta \int_s d\mathbf{r} \, \mathbf{f}_i(\mathbf{r}) \cdot \int _s \overline{\mathbf{G}}(\mathbf{r},\mathbf{r}') \cdot \mathbf{f}_j(\mathbf{r}') \, d\mathbf{r}'
\end{aligned}
\end{equation}

As suggested by Harrington and Mautz \cite{Harrington1971a}, the CMs of arbitrary PEC objects can be obtained by solving the GEP give by
\begin{equation}\label{GEP}
\bar{\mathbf{X}}  \mathbf{J}_n = \lambda_n \bar{\mathbf{R}}   \mathbf{J}_n 
\end{equation}
where $\bar{\mathbf{R}}$ and $\bar{\mathbf{X}}$ are the resistance and reactance matrices, respectively. They are the real and imaginary parts of $\bar{\mathbf{Z}} $, i.e.,
\begin{equation}
\bar{\mathbf{Z}}  =\bar{\mathbf{R}}  + i\bar{\mathbf{X}} 
\end{equation}
The eigenpairs $\lambda_n$ and $\mathbf{J}_n $ correspond to the characteristic values and characteristic currents, respectively. The characteristic value $\lambda_n$ is important as $\vert \lambda_n \vert$ indicates the modal behavior. When $\lambda_n=0$, the corresponding $\mathbf{J}_n $ is at resonance which is efficient in radiating energy. When $\vert \lambda_n \vert$ is large, $\mathbf{J}_n$ is an inefficient radiating mode which stores energy in the object. More explicitly, when $\lambda_n>0\, (\lambda_n<0)$, $\mathbf{J}_n$ stores magnetic (electric) energy.
    
\section{Multilevel Fast Multipole Algorithm}
\subsection{Addition Theorem}
The standard MLFMA is based on the addition theorem. Letting $\mathbf{r}_i$ and $\mathbf{r}_j$ be the observation point and source point located in group $m$ and $m^\prime$, respectively, the spatial vector from the source point to the observation point is written as \cite{Chew2000}
\begin{equation}
\begin{aligned}
\mathbf{r}_{ij} = \mathbf{r}_i - \mathbf{r}_j &= \mathbf{r}_i - \mathbf{r}_m + \mathbf{r}_m - \mathbf{r}_{m'} + \mathbf{r}_{m'} - \mathbf{r}_j \\
& = \mathbf{r}_{im} + \mathbf{r}_{mm'} + \mathbf{r}_{m'j} 
\end{aligned}
\end{equation} 
where $\mathbf{r}_m$ and $\mathbf{r}_{m'}$ represent the centers of groups $m$ and $m'$. If groups $m$ and $m'$ are not nearby, i.e., $\vert \mathbf{r}_{im} + \mathbf{r}_{m'j} \vert < \vert \mathbf{r}_{mm'}\vert$, their interaction is considered as far interaction, and the spherical-wave function can be expanded using the addition theorem as
\begin{equation}
\frac{e^{i k r_{ij}}}{ r_{ij}} = \frac{i k}{4 \pi} \int d^2 \hat k \, e^{i \mathbf{k}\cdot (\mathbf{r}_{im} + \mathbf{r}_{m'j})}\, \alpha_{mm'}\left(\mathbf{k}, \mathbf{r}_{mm'}\right)
\end{equation}
where $r_{ij}=\vert \mathbf{r}_i - \mathbf{r}_j\vert $, and the translator $\alpha_{mm'}$ is defined as
\begin{equation}
\alpha_{mm'}\left(\mathbf{k} , \mathbf {r}_{mm'}\right) = \sum_{l=0}^L i^{l} (2l+1) h_l^{(1)} (kr_{mm'}) P_l\left(\hat k \cdot \hat r_{mm'} \right)
\end{equation}
with $L$ as the truncation number of an infinite series, $h_l^{(1)}(\cdot)$ the first kind spherical Hankel function of order $l$, and $P_l(\cdot)$ the Legendre polynomial of degree $l$.

\subsection{Two-Level Algorithm}
In the two-level algorithm, the matrix-vector product is calculated as \cite{Chew2000}
\begin{equation}\label{ZI}
\begin{aligned}
\sum_{j=1}^N Z_{ij} I_j = & \sum_{m'\in B_{m}} \sum_{j\in G_{m'}} Z_{ij}I_j \\
& + \frac{ik}{4\pi}\int d^2 \hat k \, \mathbf{V}_{fim}(\hat k)  \cdot \sum_{m' \notin B_m}\alpha_{mm'}(\mathbf{k}, \mathbf{r}_{mm'}) \\ &\quad \cdot \sum_{j\in G_{m'}} \mathbf{V}_{sm'j}(\hat k) I_j,   \qquad i\in G_m
\end{aligned} 
\end{equation}
where
\begin{equation}
\mathbf{V}_{fim}(\hat k) = \int_S d\mathbf{r} \left(\bar{\mathbf{I}} - \hat k \hat k \right)\cdot \mathbf{f}_i (\mathbf{r}_{im}) e^{i\mathbf{k}\cdot \mathbf{r}_{im}}
\end{equation}
\begin{equation}
\mathbf{V}_{sm'j}(\hat k) = \int_S d\mathbf{r} \left(\bar{\mathbf{I}} - \hat k \hat k \right)\cdot \mathbf{f}_j (\mathbf{r}_{m'j}) e^{i\mathbf{k}\cdot \mathbf{r}_{m'j}}
\end{equation}
and $G_m$ denotes all elements in group $m$, $B_m$ denotes all nearby groups of group $m$ (including itself). The first term corresponds to the contribution of near interaction where the impedance submatrix is explicitly generated, and the second term corresponds to the contribution of far interaction evaluated by the fast multipole method (FMM). Moreover, for the far interaction, the summation $\sum_{j\in G_{m'}} \mathbf{V}_{sm'j} I_j$ represents the aggregation of all unknowns in group $m'$ to the group center, the summation $\sum_{m'}\alpha_{mm'}$ represents the translation from group $m'$ to $m$, and $\mathbf{V}_{fim}\cdots$ represents the disaggregation from the center of group $m$ to all unknowns in this group.   

\subsection{Multilevel Algorithm}
The two-level FMM is extended to a multilevel algorithm, i.e., MLFMA, to further reduce the cost of matrix-vector product operations. In MLFMA, the object $S$ is first enclosed by a box that contains $S$, which is referred to as level-$0$. This box is then partitioned into eight smaller ones to form level-$1$. Each level-$1$ box is again partitioned into eight smaller ones to form level-$2$,  and recursively continued until the finest level $L_f$ with a box size of $0.2\lambda$ to $0.5\lambda$ is formed. This procedure results in an oct-tree structure, where level-$2$ is taken as the coarsest level for best efficiency. The far interaction at the coarsest level is implemented using MLFMA where interpolation and anterpolation are adopted in the aggregation from the finest level to the coarsest level and  the disaggregation from the coarsest level to the finest level, respectively. This is because the radiation patterns become increasingly richer as one progresses from the finest level to the coarsest level. The radiation patterns are sampled to perform the integration over the Ewald sphere in (\ref{ZI}). A well-documented approach is to apply the uniform sampling and Gauss-Legendre sampling in the $\phi$ and $\theta$ coordinates, respectively, and local Langrange interpolation between levels \cite{Chew2000}.

\section{MLFMA for CM Analysis}
For large scale CM analysis, the desired eigenpairs of (\ref{GEP}) are normally solved iteratively. The ARnoldi PACKage (ARPACK), based on the implicitly restarted Arnoldi (IRA) method or the Lanczos variant for symmetric matrices, is appropriate for calculating a few eigenpairs of large sparse or structured matrices, where a matrix-vector product requires only $O(N)$ floating point operations \cite{Sorensen}. However, when ARPACK is applied to the dense matrix pencil $\left(\bar{\mathbf{X}},\bar{\mathbf{R}}\right)$ in (\ref{GEP}), the matrix-vector product is challenged by the $O(N^2)$ complexity and storage requirement. 

In most applications, the eigenvalues of great interest are those with small magnitudes which correspond to efficient radiating modes. Since it is easy for the IRA method to converge to large eigenvalues, the desired eigenpairs can be found by solving the following standard eigenvalue problem
\begin{equation}\label{SEP}
\bar{\mathbf{X}}^{-1}  \bar{\mathbf{R}}  \mathbf{J}_n= \frac{1}{\lambda_n} \mathbf{J}_n
\end{equation}
The Arnoldi process for (\ref{SEP}) calls for fast matrix-vector products in the form of $\bar{\mathbf{X}}^{-1}\bar{\mathbf{R}}\mathbf{u}$ with an arbitrary vector $\mathbf{u}$ to construct Arnoldi vectors. Even though matrices $\bar{\mathbf{R}}$ and $\bar{\mathbf{X}}$ can be explicitly generated if tremendous memory is available, the matrix vector product $\bar{\mathbf{X}}^{-1}\mathbf{v}$ where $\mathbf{v}=\bar{\mathbf{R}}  \mathbf{u}$ can only be iteratively computed, as the complexity of inverting $\bar{\mathbf{X}}$ is $O(N^3)$.  

It is well known that MLFMA reduces the complexity and memory cost of performing $\bar{\mathbf{Z}} \mathbf{u}$ to $O(N\log N)$. 
The standard MLFMA can be easily modified to expedite the required matrix vector products.
To calculate $\bar{\mathbf{X}}  \mathbf{u} $ with an arbitrary $\mathbf{u}$, the following plane-wave decomposition is used provided $\vert \mathbf{r}_{im} + \mathbf{r}_{m'j} \vert < \vert \mathbf{r}_{mm'}\vert$ 
\begin{equation}\label{CosKer}
\frac{\cos( k  r_{ij})}{ r_{ij} } = \frac{i k}{4 \pi} \int d^2 \hat k \, e^{i \mathbf{k}\cdot (\mathbf{r}_{im} + \mathbf{r}_{m'j})}\, \alpha_{mm'}^\mathfrak{S}\left(\mathbf{k}, \mathbf{r}_{mm'}\right)
\end{equation}
where
\begin{equation}
\alpha_{mm'}^\mathfrak{S}\left(\mathbf{k}, \mathbf{r}_{mm'}\right)=i \sum_{l=0}^L i^{l}(2l+1)  y_l (kr_{mm'})  P_l\left(\hat k \cdot \hat r_{mm'} \right)
\end{equation}
with the spherical Neumann function  $y_l(\cdot)$ of order $l$. 
Furthermore, the reactance matrix $\bar{\mathbf{X}} $ is given by
\begin{equation}
X_{ij} = \Im m \lbrack Z_{ij} \rbrack = k\eta \left\langle \mathbf{f}_i, \left\lbrack\overline{\mathbf{I}}-\frac{\nabla\nabla'}{k^2}\right\rbrack \frac{\cos\left(kR \right)}{R}, \mathbf{f}_j\right\rangle_S 
\end{equation} 
The required matrix-vector product in terms of a two-level algorithm is  obtained as
\begin{equation}\label{XI}
\begin{aligned}
\sum_{j=1}^N X_{ij} u_j = & \sum_{m'\in B_{m}} \sum_{j\in G_{m'}} X_{ij} u_j \\
& + \frac{ik}{4\pi}\int d^2 \hat k \, \mathbf{V}_{fim}(\hat k)  \cdot \sum_{m' \notin B_m}\alpha_{mm'}^\mathfrak{S}(\mathbf{k}, \mathbf{r}_{mm'}) \\ &\quad \cdot \sum_{j\in G_{m'}} \mathbf{V}_{sm'j}(\hat k) u_j,   \qquad i\in G_m
\end{aligned} 
\end{equation}
where $G_m$, $B_m$, $\mathbf{V}_{fim}$ and $\mathbf{V}_{sm'j}$ are defined the same as before. The corresponding multilevel algorithm simply follows the same scheme used in the standard MLFMA. 

The calculation of $\bar{\mathbf{R}}  \mathbf{u}$ is less involved as there is no violation of the following addition theorem  
\begin{equation}\label{SinKer}
\frac{\sin( k  r_{ij})}{ r_{ij} } = \frac{i k}{4 \pi} \int d^2 \hat k \, e^{i \mathbf{k}\cdot (\mathbf{r}_{im} + \mathbf{r}_{m'j})}\, \alpha_{mm'}^\mathfrak{R}\left(\mathbf{k}, \mathbf{r}_{mm'}\right)
\end{equation}
where
\begin{equation}
\alpha_{mm'}^\mathfrak{R}\left(\mathbf{k}, \mathbf{r}_{mm'}\right)=-i \sum_{l=0}^L i^{l}(2l+1) j_l  (kr_{mm'})  P_l\left(\hat k \cdot \hat r_{mm'} \right)
\end{equation}
with the spherical Bessel function  $j_l(\cdot)$ of order $l$. Moreover, the resistance matrix $\bar{\mathbf{R}}$ is given by  
\begin{equation}
R_{ij} = \Re e\lbrack Z_{ij} \rbrack =  -k\eta \left\langle \mathbf{f}_i, \left\lbrack\overline{\mathbf{I}}-\frac{\nabla\nabla'}{k^2}\right\rbrack \frac{\sin\left(kR\right)}{R}, \mathbf{f}_j\right\rangle_S \\
\end{equation} 
The required matrix-vector product in terms of a two-level algorithm is obtained as
\begin{equation}\label{RI}
\begin{aligned}
\sum_{j=1}^N R_{ij} u_j = &  \sum_{j\in G_{m}} R_{ij} u_j \\
& + \frac{ik}{4\pi}\int d^2 \hat k \, \mathbf{V}_{fim}(\hat k)  \cdot \sum_{m' \neq m}\alpha_{mm'}^\mathfrak{R}(\mathbf{k}, \mathbf{r}_{mm'}) \\ &\quad \cdot \sum_{j\in G_{m'}} \mathbf{V}_{sm'j}(\hat k) u_j,   \qquad i\in G_m
\end{aligned} 
\end{equation}
where the FMM calculation involves both the nearby and nonnearby groups of group $m$ (except for itself). With regard to the multilevel algorithm, the standard MLFMA bookkeeping can still be employed for the least modification of the existing code. However, a more efficient implementation conducts translation between nearby groups which belong to the same parent group at each level. Besides, the coarsest level is taken to be level-$1$ instead of level-$2$.

When the problem scale is large, the reactance matrix $\bar{\mathbf{X}}$ becomes ill-conditioned. It may take many iterations for the generalized minimal residual (GMRES) method or biconjugate gradient stabilized (BiCGSTAB) method to converge to a prescribed error tolerance as $\bar{\mathbf{X}}^{-1}  \mathbf{u}$ is solved. On the other hand, efficient preconditioning schemes for solving $\bar{\mathbf{X}} ^{-1}\mathbf{u}$ are not well documented. An alternative is to consider the equivalent  eigenvalue problem of (\ref{GEP}) given by
\begin{equation}\label{GEP1}
\bar{\mathbf{Z}}  \mathbf{J}_n  = (1+i\lambda_n) \bar{\mathbf{R}}  \mathbf{J}_n 
\end{equation}
Since $\lambda_n$ are real numbers in the CM theory, $\lambda_n$  close to zero can be easily found with IRA by solving eigenvalues with large magnitude of the standard eigenvalue problem
\begin{equation}\label{SEP1}
\bar{\mathbf{Z}}^{-1} \bar{\mathbf{R}}  \mathbf{J}_n = \frac{1}{1+i\lambda_n} \mathbf{J}_n 
\end{equation}
Although IRA may require more precision in solving $\bar{\mathbf{Z}} ^{-1}  \mathbf{u}$ than $\bar{\mathbf{X}}^{-1}  \mathbf{u}$ to yield accurate eigenpairs, (\ref{SEP1}) is still favored in this study as there exist efficient preconditioners for accelerating the EFIE solutions  $\bar{\mathbf{Z}}^{-1} \mathbf{u}$. These preconditioners may not work well in preconditioning $\bar{\mathbf{X}}^{-1}  \mathbf{u}$ as the Green's function $\cos(kR)/R$ is noncausal. An SAI preconditioner is implemented in this study following \cite{Jeonghwa} due to its simplicity. In MLFMA, the impedance submatrix $\bar{\mathbf{Z}}_{near} $ for near interaction is explicitly evaluated and stored. The sparse preconditioning matrix $\bar{\mathbf{P}}$ with a prescribed sparsity pattern is constructed by minimizing  $ \left\|\bar{\mathbf{I}} - \bar{\mathbf{Z}}_{near} \bar{\mathbf{P}}  \right \|_F$ where the subscript $F$ denotes the Frobenius norm.

A more involved but even better preconditioner is the CMP that projects $\bar{\mathbf{Z}}   \mathbf{x} =  \mathbf{u}$ to form \cite{Andriulli}
\begin{equation}
\widetilde{\bar{\mathbf{Z}}}  \bar{\mathbf{G}}^{-1}  \bar{\mathbf{Z}}  \mathbf{x}  = \widetilde{\bar{\mathbf{Z}}}  \bar{\mathbf{G}}^{-1}  \mathbf{u}
\end{equation} 
where 
\begin{equation}
\widetilde{Z}_{ij} =ik\eta \left \langle  \widetilde{\mathbf{f}}_i, \bar{\mathbf{G}}, \widetilde{\mathbf{f}}_j \right \rangle_S
\end{equation}
and the Gram matrix
\begin{equation} 
G_{ij} =\left \langle \hat n \times \mathbf{f}_i,  \widetilde{\mathbf{f}}_j \right \rangle_S
\end{equation}
where $\widetilde{\mathbf{f}}_i$ and $\widetilde{\mathbf{f}}_j$ are chosen as Buffa-Christiansen (BC) basis functions. Due to the Calderon identity, $\widetilde{\bar{\mathbf{Z}}}  \bar{\mathbf{G}} ^{-1}  \bar{\mathbf{Z}} $ is well-conditioned for a uniform mesh. When the mesh is nonuniform, diagonal preconditioners are applied to guarantee low condition numbers for both $\bar{\mathbf{G}} $ and $\widetilde{\bar{\mathbf{Z}}}  \bar{\mathbf{G}} ^{-1}  \bar{\mathbf{Z}}$. The reason for keeping $\bar{\mathbf{G}} $ well-conditioned is that one needs to compute $\bar{\mathbf{G}}^{-1}  \mathbf{u}$ iteratively provided $\bar{\mathbf{G}}$ is sparse.  Moreover, MLFMA has been embedded within CMP, which offers promise for solving large scale EFIEs with good accuracy using only a few iterations \cite{Peeters}.

The IRA method is less appropriate for solving (\ref{SEP}) due to the lack of efficiency in preconditioning $\bar{\mathbf{X}} ^{-1}  \mathbf{u}$. However, the GEP (\ref{GEP}) can be solved with the JDQZ method \cite{Fokkema}. Without converting (\ref{GEP}) to (\ref{SEP}), the correction vector $\mathbf{t}$ is obtained by approximately solving the Jacobi-Davidson correction equation
\begin{equation}\label{JDCE}
\begin{split}
&\left(\bar{\mathbf{I}} - \bar{\mathbf{R}}  \mathbf{u}  \mathbf{u}^*\right)  \left(\bar{\mathbf{X}} - \sigma \bar{\mathbf{R}} \right)  \left( \bar{\mathbf{I}}  - \mathbf{u} \mathbf{u}^* \bar{\mathbf{R}} \right) \mathbf{t} \\
& \hspace{1in}  = - \left( \bar{\mathbf{X}} - \sigma \bar{\mathbf{R}} \right)\mathbf{u}, \quad \mathbf{t}^* \bar{\mathbf{R}}  \mathbf{u} =0 
\end{split}
\end{equation}
where $\left(\sigma, \mathbf{u} = \bar{\mathbf{V}} \mathbf{s}\right)$ is a solution of 
\begin{equation}
\bar{\mathbf{V}}^* \bar{\mathbf{X}} \bar{\mathbf{V}} \mathbf{s} = \sigma \bar{\mathbf{V}}^* \bar{\mathbf{R}} \bar{\mathbf{V}} \mathbf{s} = \sigma \mathbf{s}
\end{equation}
with the search space $\bar{\mathbf{V}}$ containing several $R$-orthogonal basis. The key is that even when (\ref{JDCE}) is solved with low accuracy, JDQZ still gives rise to correct solutions of (\ref{GEP}). However, JDQZ becomes fairly slow in this case. Hence, it is also beneficial if efficient preconditioners are designed to accelerate the solutions of (\ref{JDCE}).

\section{Numerical Results} 
\begin{figure}[!t]
\centering
\includegraphics[scale=0.45]{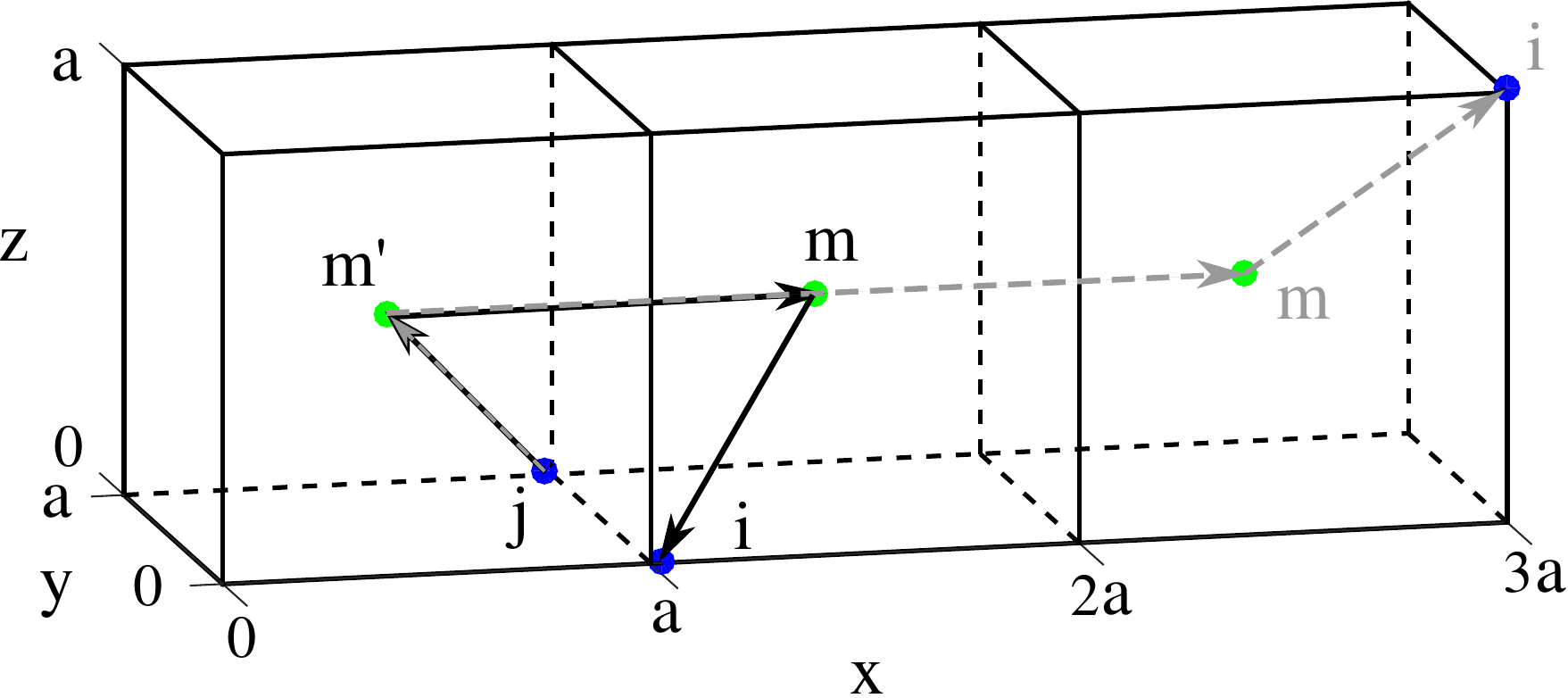}
\caption{Point positions.} \label{box}
\end{figure}

\begin{figure}[!t]
\centering
\includegraphics[scale=0.4]{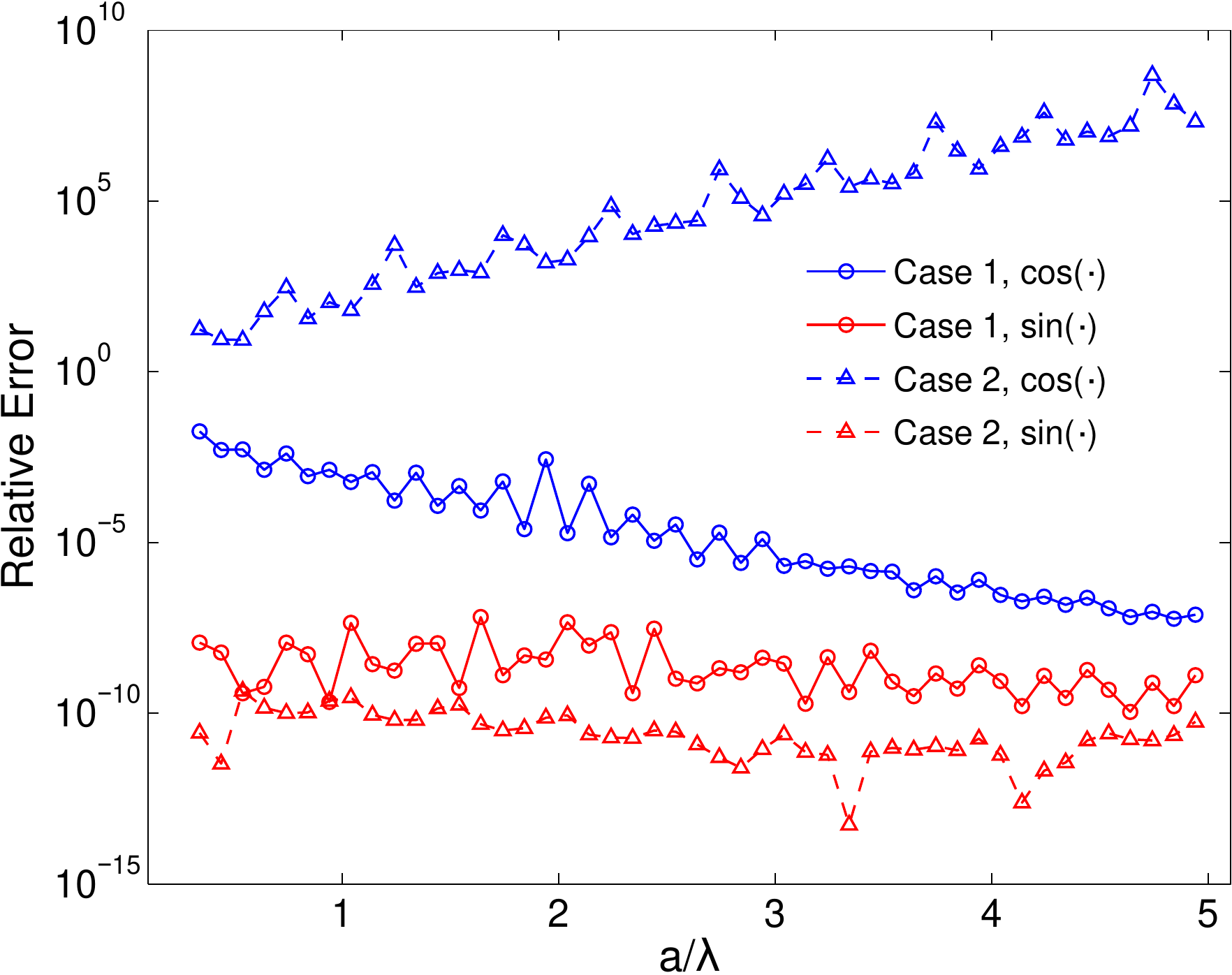}
\caption{Relative error of plane wave decomposition for the real and imaginary parts of Green's function with respect to $a/\lambda$.} \label{pointTestError}
\end{figure}

Numerical tests are first performed to validate Equations (\ref{CosKer}) and (\ref{SinKer}). Two cases are considered as shown in Fig.~\ref{box}, where the box size is $a$. In both cases, the source point $j$ and source group center $m'$ are located at $(0.9999a,\thinspace 1.0a,\thinspace 0.0)$ and $(0.5a,\thinspace 0.5a,\thinspace 0.5a)$, respectively. The observation point $i$ and its group center $m$ in Case 1 (dashed arrows) are located at $(2.9999a,\thinspace 0.0,\thinspace 1.0a)$ and $(2.5a,\thinspace 0.5a,\thinspace 0.5a)$, respectively, such that $\vert \mathbf{r}_{im} + \mathbf{r}_{m'j} \vert < \vert \mathbf{r}_{mm'}\vert$. The relative errors between left and right hand sides of (\ref{CosKer}) and  (\ref{SinKer}) are calculated with respect to  $a/\lambda$, where $\lambda$ is the wavelength. The truncation number $L$ is chosen to be 
\begin{equation}
L\approx kd + 1.8 d_0^{2/3} (kd)^{1/3}
\end{equation}
where the number of digits of accuracy $d_0$ is set to $3$, and $d=\sqrt{3} a$ in this test.
As shown in Fig.~\ref{pointTestError}, excellent agreement is observed except for (\ref{CosKer}) or $\cos(\cdot)$ at lower frequencies ($a<0.5\lambda$), which can be resolved with the low-frequency fast multipole algorithm (LF-FMA) \cite{Zhao00}. We then place the observation point $i$ and its group center $m$  at $(1.0001a,\thinspace 0.0,\thinspace 0.0)$ and $(1.5a,\thinspace 0.5a,\thinspace 0.5a)$, respectively, in Case 2 (solid arrows) such that $\vert \mathbf{r}_{im} + \mathbf{r}_{m'j} \vert > \vert \mathbf{r}_{mm'}\vert$. As expected, (\ref{CosKer}) is no longer valid due to the violation of addition theorem. However,  the regular part of the Green's function, i.e., $\sin(\cdot)$ can be computed with excellent accuracy  [Fig.~\ref{pointTestError}]. This holds even when $i$ and $j$ are very close to each other. Although the results computed by $(\ref{CosKer})$ and $(\ref{SinKer})$ are not purely real, the imaginary parts are negligible as they are normally $10^{10}$ times smaller than the real parts.

\begin{figure}[!t]
\centering
\subfigure[]{\includegraphics[scale=0.21]{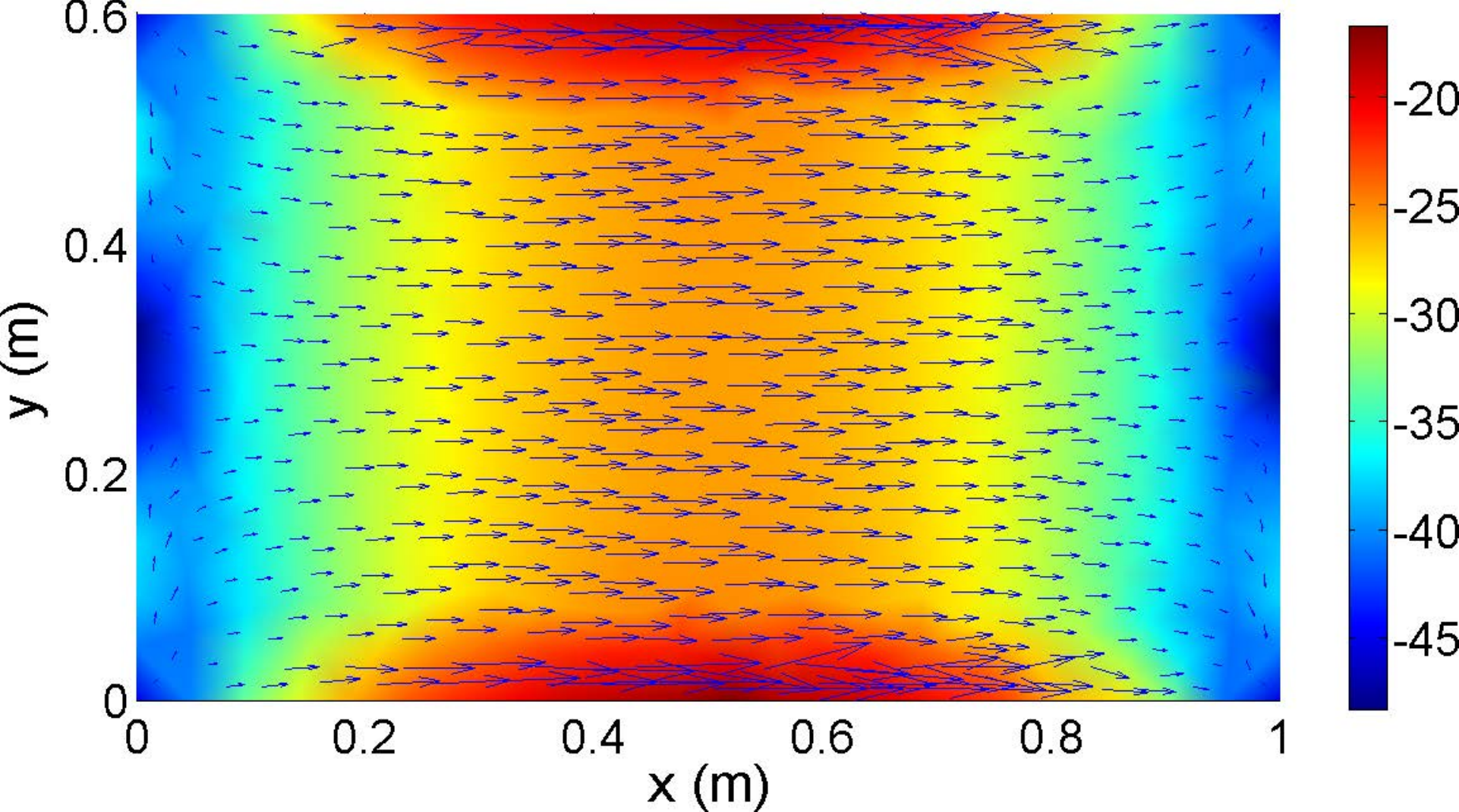} \label{J1}}
\subfigure[]{\includegraphics[scale=0.21]{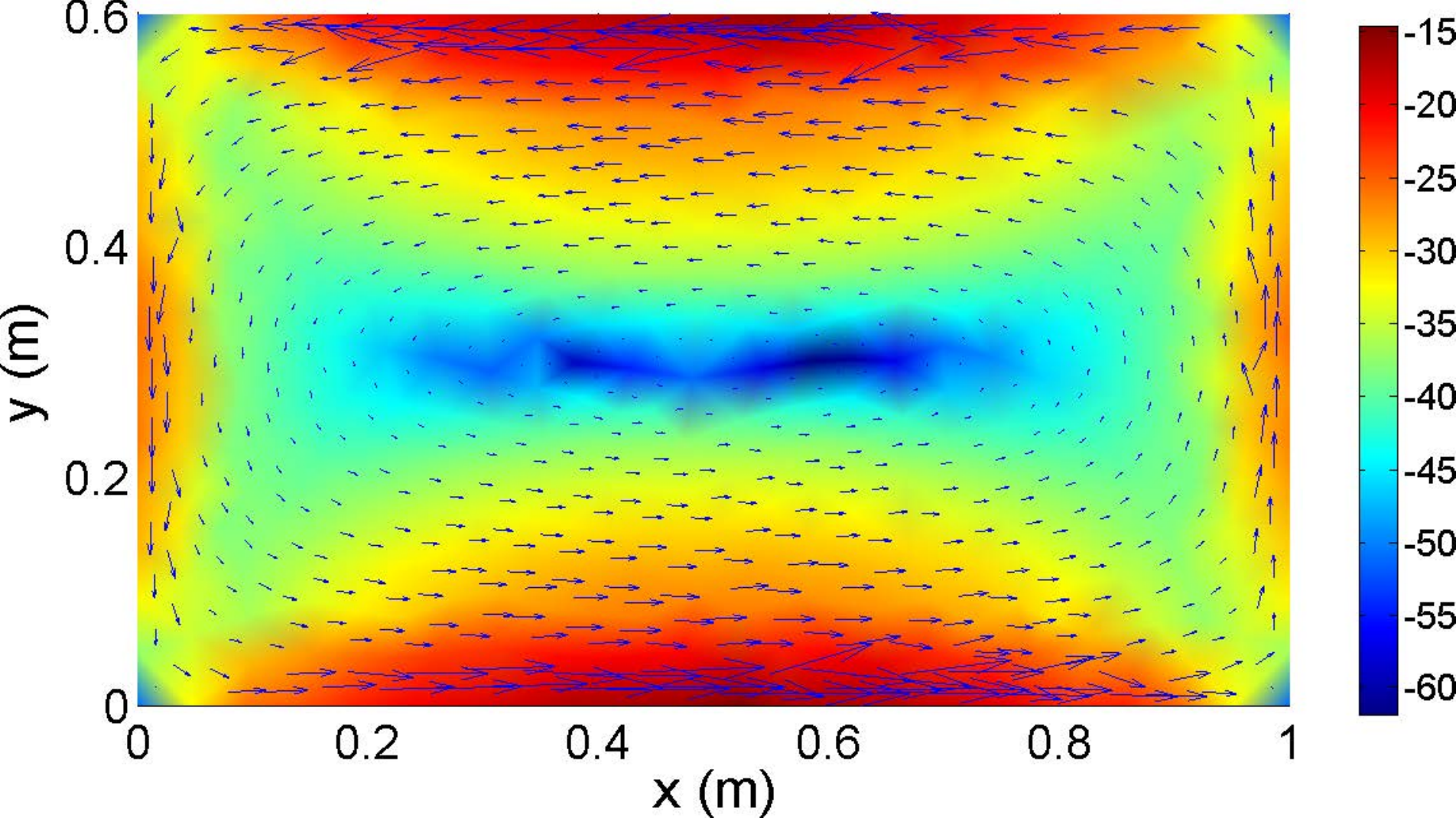} \label{J2}} 
\subfigure[]{\includegraphics[scale=0.21]{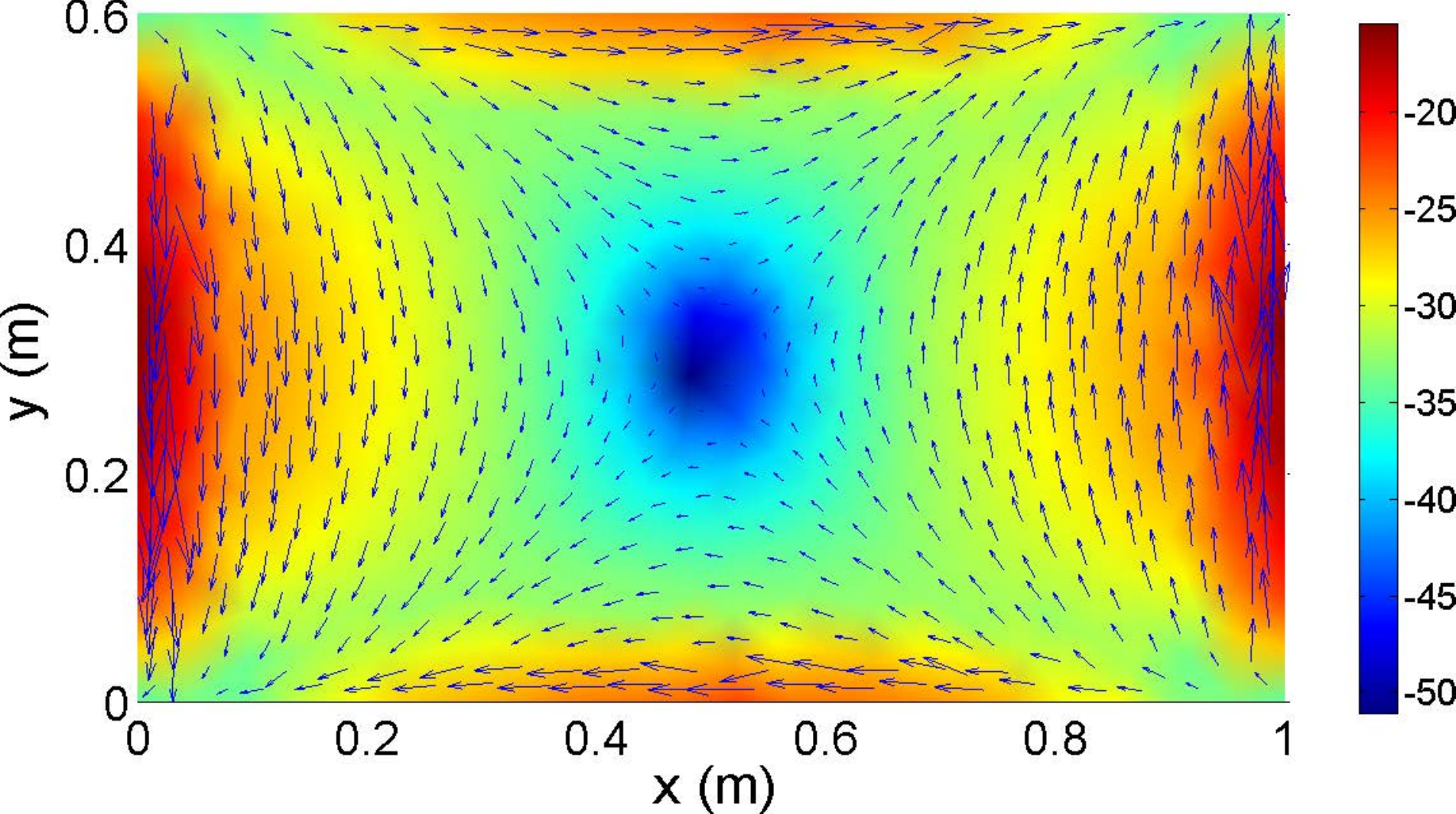} \label{J3}} 
\subfigure[]{\includegraphics[scale=0.21]{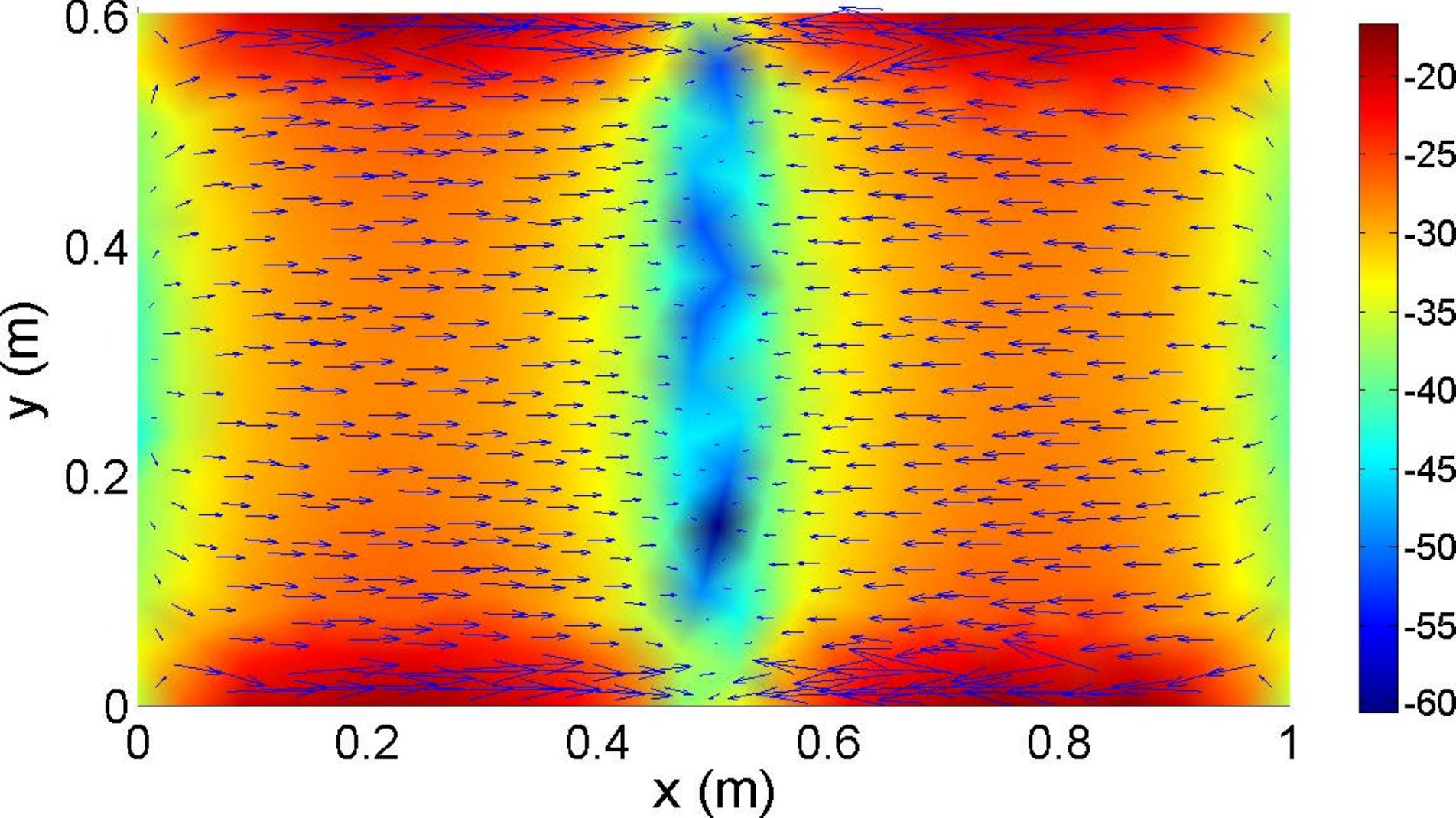} \label{J4}} 
\caption{Current patterns of plate characteristic modes: (a) $\mathbf{J}_1$. (b) $\mathbf{J}_2$. (c) $\mathbf{J}_3$. (d) $\mathbf{J}_4$} \label{CM Plate}
\end{figure}

In this study, iterations in IRA or JDQZ (GMRES or BiCGSTAB) are referred to as outer (inner) iterations. To further test the validity of the proposed method, we perform the CM analysis to a rectangular PEC plate. The dimension of the plate is $1.0$~m $\times$ $0.6$~m, and the frequency varies from $220$~MHz to $380$~MHz at a step of $20$~MHz. Totally $947$ unknowns are generated in this simulation, rendering MLFMA a three-level algorithm. By solving (\ref{SEP1}) with MLFMA and IRA, four CMs whose current patterns at $300$~MHz are plotted in Fig.~\ref{CM Plate}. IRA performs $20$ outer iterations to obtain $5$ desired modes, where GMRES is used to compute $\bar{\mathbf{Z}}^{-1} \mathbf{u}$ with a tolerance of $10^{-3}$. The threshold parameters $\epsilon_1$, $\epsilon_2$ and $\epsilon_3$  in the SAI preconditioner \cite{Jeonghwa} are specified to be $0.01$, $0.014$ and $0.18$, respectively. 

\begin{figure}[!t]
\centering
\includegraphics[scale=0.45]{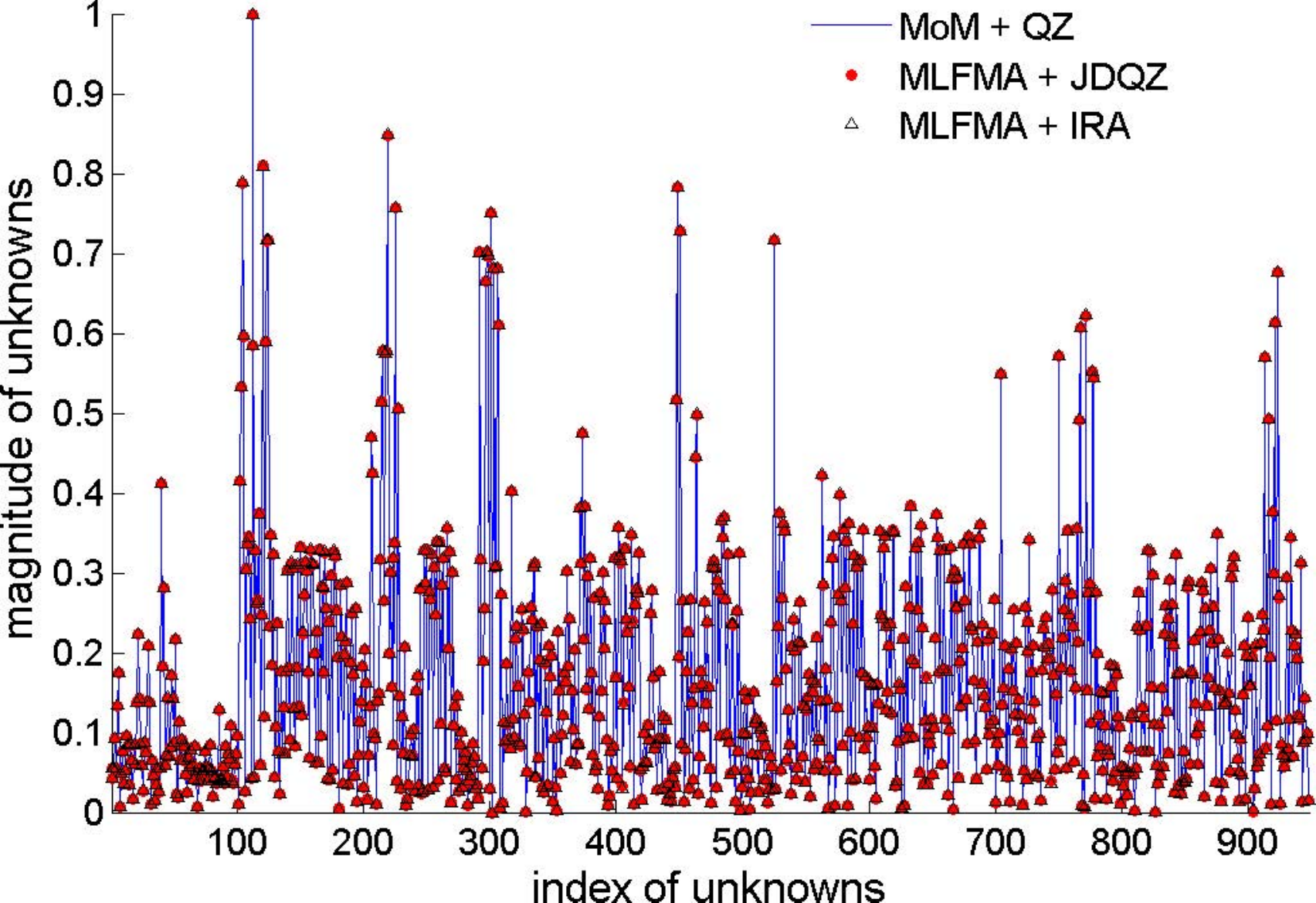}
\caption{Magnitude of unknowns in current mode $\mathbf{J}_1$} \label{cmpUnk}
\end{figure}

The same current patterns are also  obtained by solving (\ref{GEP}) with MLFMA and JDQZ, or by applying the generalized Schur decomposition (also known as QZ decomposition) to (\ref{GEP}) with $\bar{\mathbf{X}}$ and $\bar{\mathbf{R}}$ generated by MoM. 
The magnitudes of $947$ unknowns in $ \mathbf{J}_1$ computed by these approaches agree well with each other, as shown in Fig.~\ref{cmpUnk}. Taking the result obtained by MoM + QZ as a reference, the magnitude error of the entire unknowns is $0.67\%$ when $\mathbf{J}_1$ is computed by MLFMA + IRA, and $0.64\%$ when computed by MLFMA + JDQZ. For the other three current modes, all relative errors are less than $1\%$.  When the tolerance for GMRES in IRA is set to $10^{-2}$, the errors are around $3\%$ for all four modes.

\begin{figure}[!t]
\centering
\subfigure[]{\includegraphics[scale=0.23]{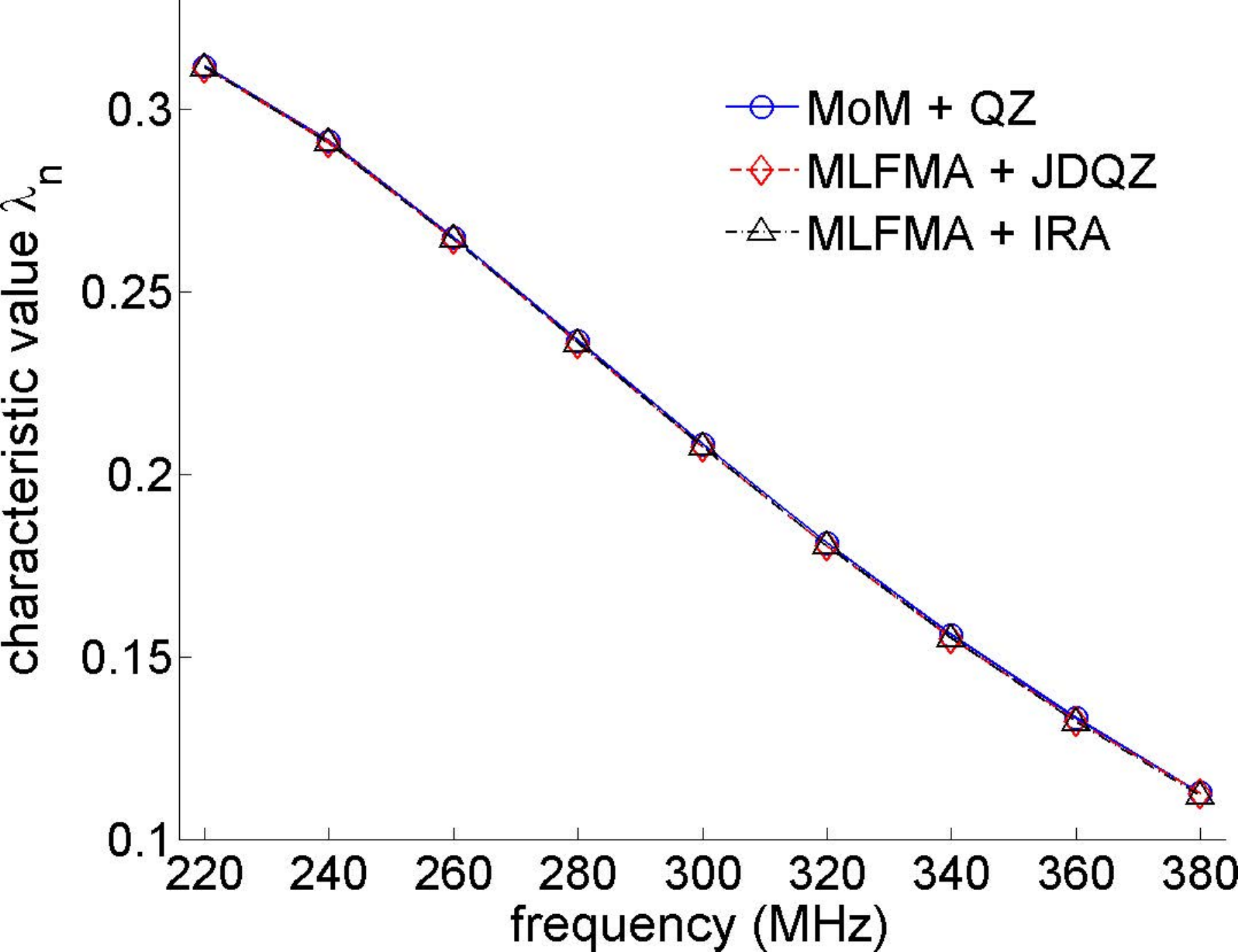} \label{CV1}}
\subfigure[]{\includegraphics[scale=0.23]{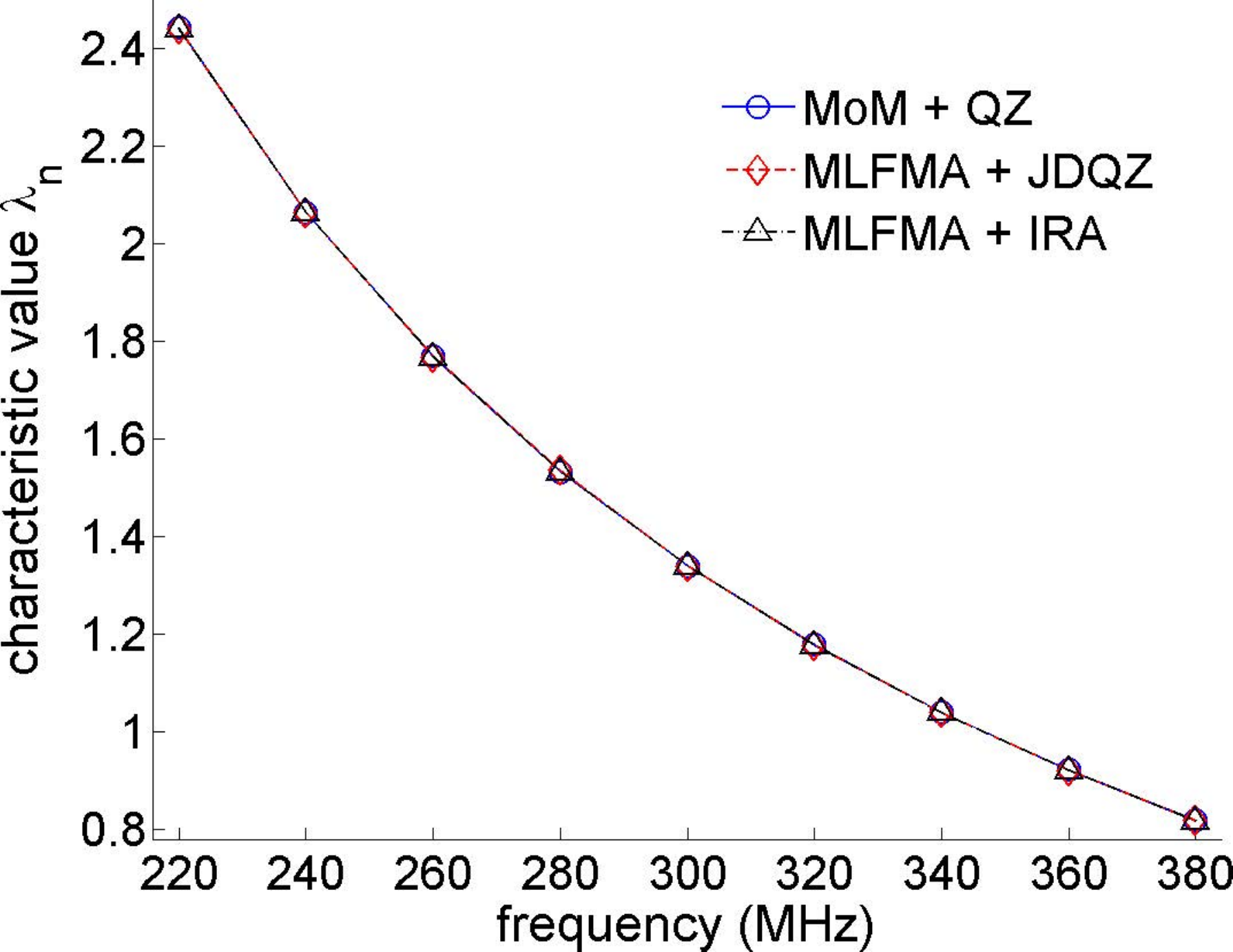} \label{CV2}} 
\subfigure[]{\includegraphics[scale=0.23]{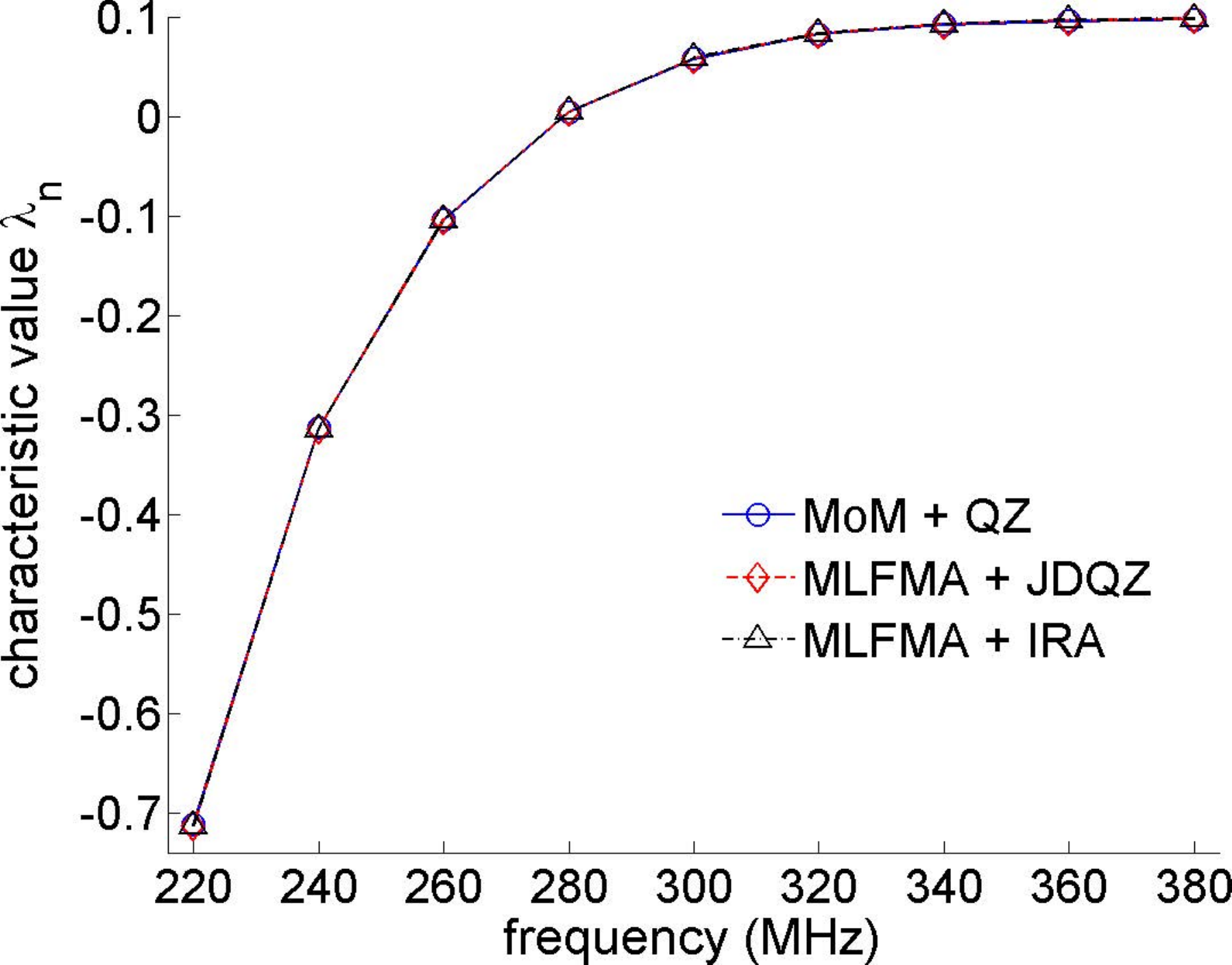} \label{CV3}} 
\subfigure[]{\includegraphics[scale=0.23]{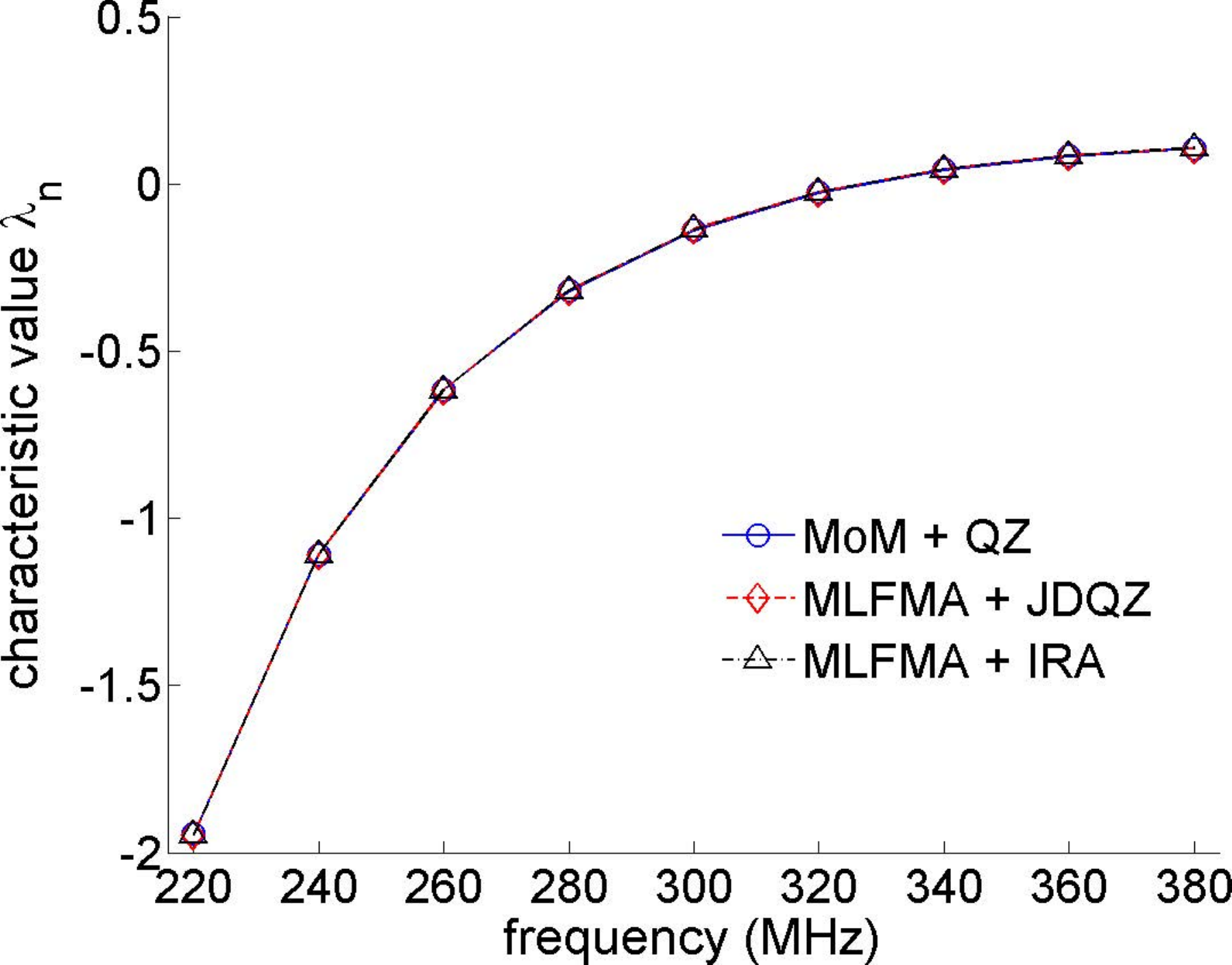} \label{CV4}} 
\caption{ Characteristic values at different frequencies: (a) $\mathbf{J}_1$. (b) $\mathbf{J}_2$. (c) $\mathbf{J}_3$. (d) $\mathbf{J}_4$} \label{CV Plate}
\end{figure}

\begin{figure}[!t]
\centering
\includegraphics[scale=0.4]{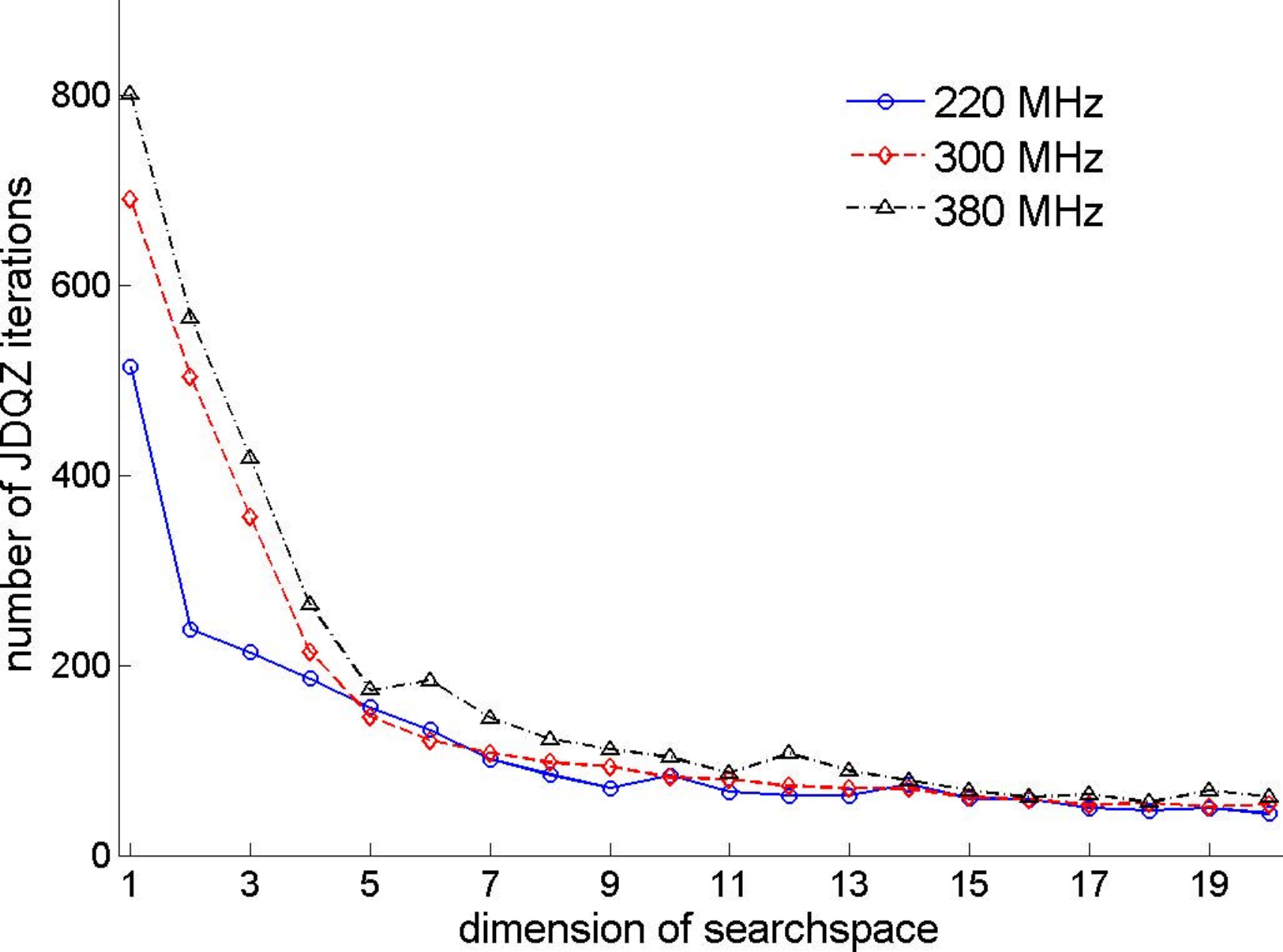}
\caption{Number of JDQZ iterations with respect to dimension of GMRES searchspace.} \label{iNjdqz}
\end{figure}

Characteristic values of $\mathbf{J}_1$ to $\mathbf{J}_4$ at different frequencies are calculated with the above approaches, where good agreement is observed from Fig.~\ref{CV Plate}. In JDQZ, one can obtain correct characteristic values even with a small number of inner iterations (dimension of searchspace in GMRES), which indicates that (\ref{JDCE}) is only solved with poor accuracy each time. Fig.~\ref{iNjdqz} shows the required number of JDQZ iterations for certain convergence tolerance, i.e., $10^{-4}$, with respect to the dimension of GMRES searchspace. Although a small searchspace reduces the computational time of each outer iteration, the tradeoff is that it increases the total number of outer iterations.

\begin{figure}[!t]
\centering
\subfigure[]{\includegraphics[scale=0.33]{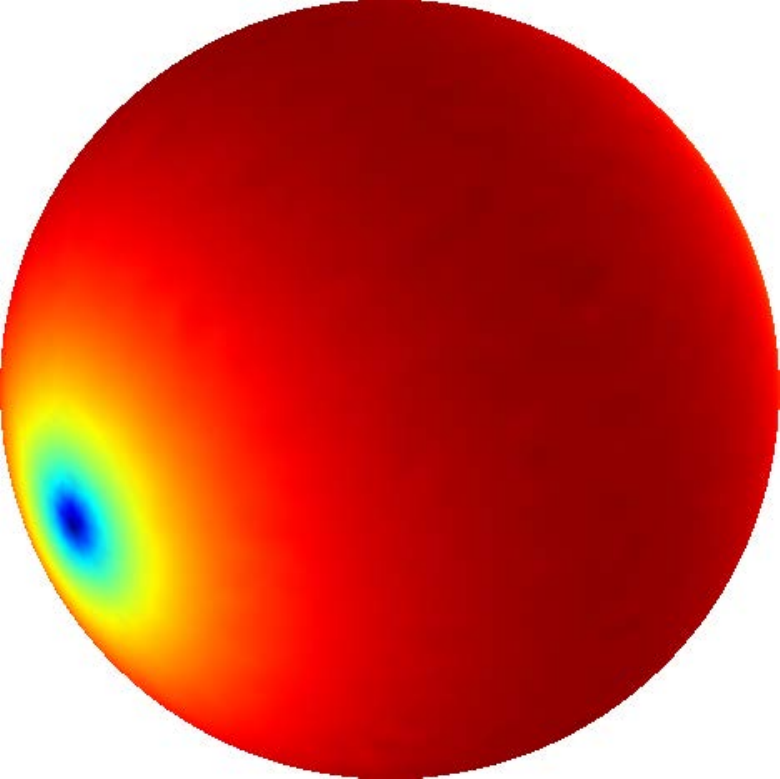} \label{CM1Sph}}
\subfigure[]{\includegraphics[scale=0.33]{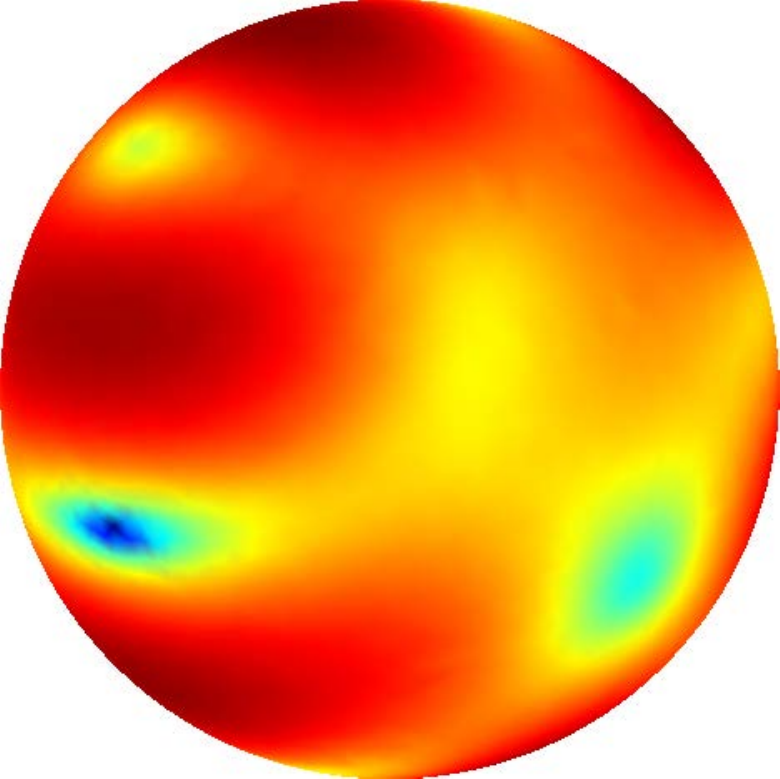} \label{CM4Sph}} 
\subfigure[]{\includegraphics[scale=0.33]{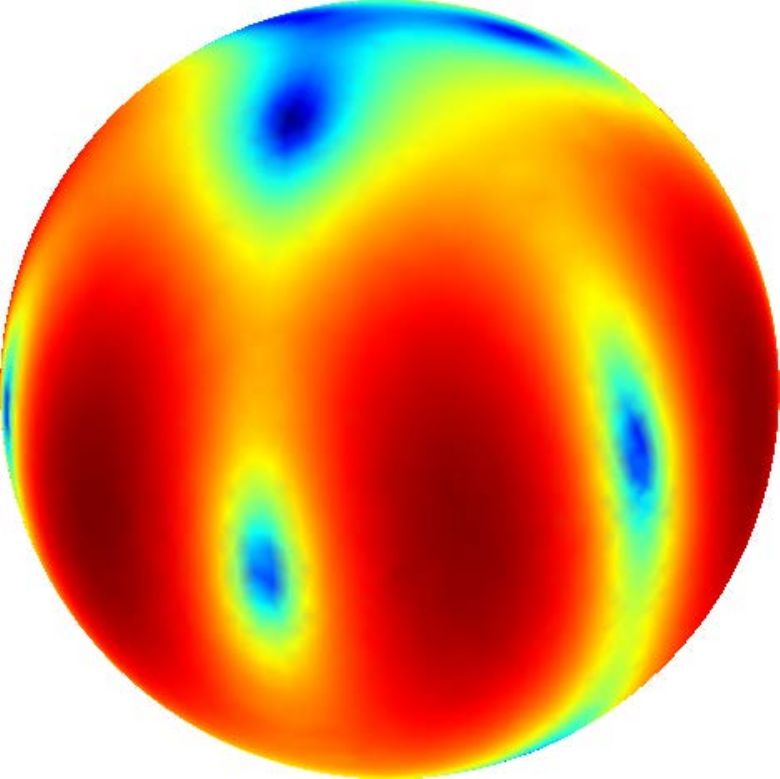} \label{CM7Sph}} 
\caption{Sphere  characteristic modes computed using MLFMA + IRA ($15\thinspace015$ unknowns): (a) $\mathbf{J}_1$. (b) $\mathbf{J}_4$. (c) $\mathbf{J}_7$.} \label{CM Sphere}
\end{figure}

We next consider a PEC sphere with a radius of $0.75\lambda$. We increase the mesh density such that  $942$, $3\thinspace 519$ and $15\thinspace 015$ unknowns are generated. A few degenerate CMs are computed with MoM + QZ and MLFMA + IRA, respectively, where the corresponding characteristic values are listed in Table \ref{TCV}. The GMRES tolerance is set to $10^{-3}$ in the IRA method. Good agreement is observed between results computed with the two schemes. In \cite{Ludick} where a few CMs of a vehicle with $9\thinspace 706$ unknowns are solved for, the memory usage is $3.54$~GB and $2.8$~GB for MoM + QZ and MoM + IRA, respectively. Hence, the MoM + QZ scheme has not been applied to the problem with $15\thinspace015$ unknowns which requires more memory.  
However, such a problem can be easily handled with the proposed MLFMA + IRA scheme on a small computer.  Certain characteristic current patterns (magnitudes) are illustrated in Fig.~\ref{CM Sphere}, which are also found from the results computed by MoM + QZ with less unknowns.  
    
 \begin{table} \centering
 \caption{Characteristic Values $\lambda_n$ of a Few Computed  Modes at Different Mesh Densities}
 \label{TCV}    
\begin{tabular}{cccccc}
 \hline
 \multirow{2}{*}{Current modes} &
 \multicolumn{2}{c}{MoM + QZ} &
  \multicolumn{3}{c}{MLFMA + IRA} \\
 \cline{2-6}
  & $942$ & $3519$ &  $942$ & $3519$ & $15015$ \\
 \hline
 $\mathbf{J}_1$ & -0.1974 & -0.2129 & -0.1963 & -0.2174 & -0.2249\\
 $\mathbf{J}_2$ & -0.1960 & -0.2131 & -0.1977 & -0.2170 & -0.2246 \\
 $\mathbf{J}_3$ & -0.1966 & -0.2132 & -0.1969 & -0.2172 & -0.2245\\
 $\mathbf{J}_4$ & 0.3081 &  0.2899  & 0.3078 & 0.2894 & 0.2857\\
 $\mathbf{J}_5$ & 0.3100 &  0.2900  & 0.3100 & 0.2913 & 0.2857\\
 $\mathbf{J}_6$ & 0.3098 &  0.2901  & 0.3098 & 0.2912 & 0.2856\\
 $\mathbf{J}_7$ & 0.3093 &  0.2903  & 0.3083 & 0.2912 & 0.2855\\
 $\mathbf{J}_8$ & 0.3089 &  0.2902  & 0.3086 & 0.2904 & 0.2855\\
 \hline
 \end{tabular}
\end{table}

\begin{figure}[!t]
\centering
\subfigure[]{\includegraphics[scale=0.32]{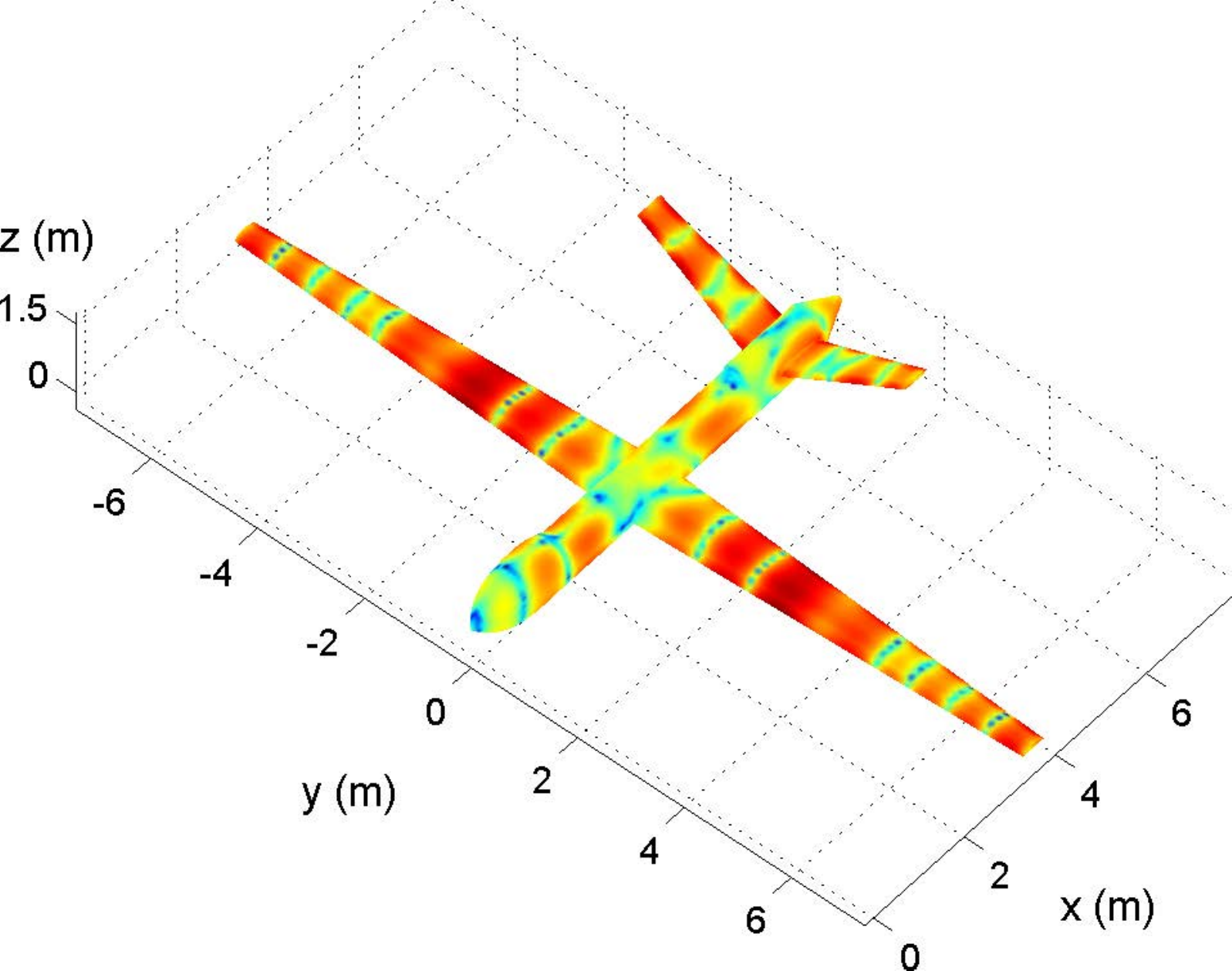} \label{CM1Pred}}
\subfigure[]{\includegraphics[scale=0.32]{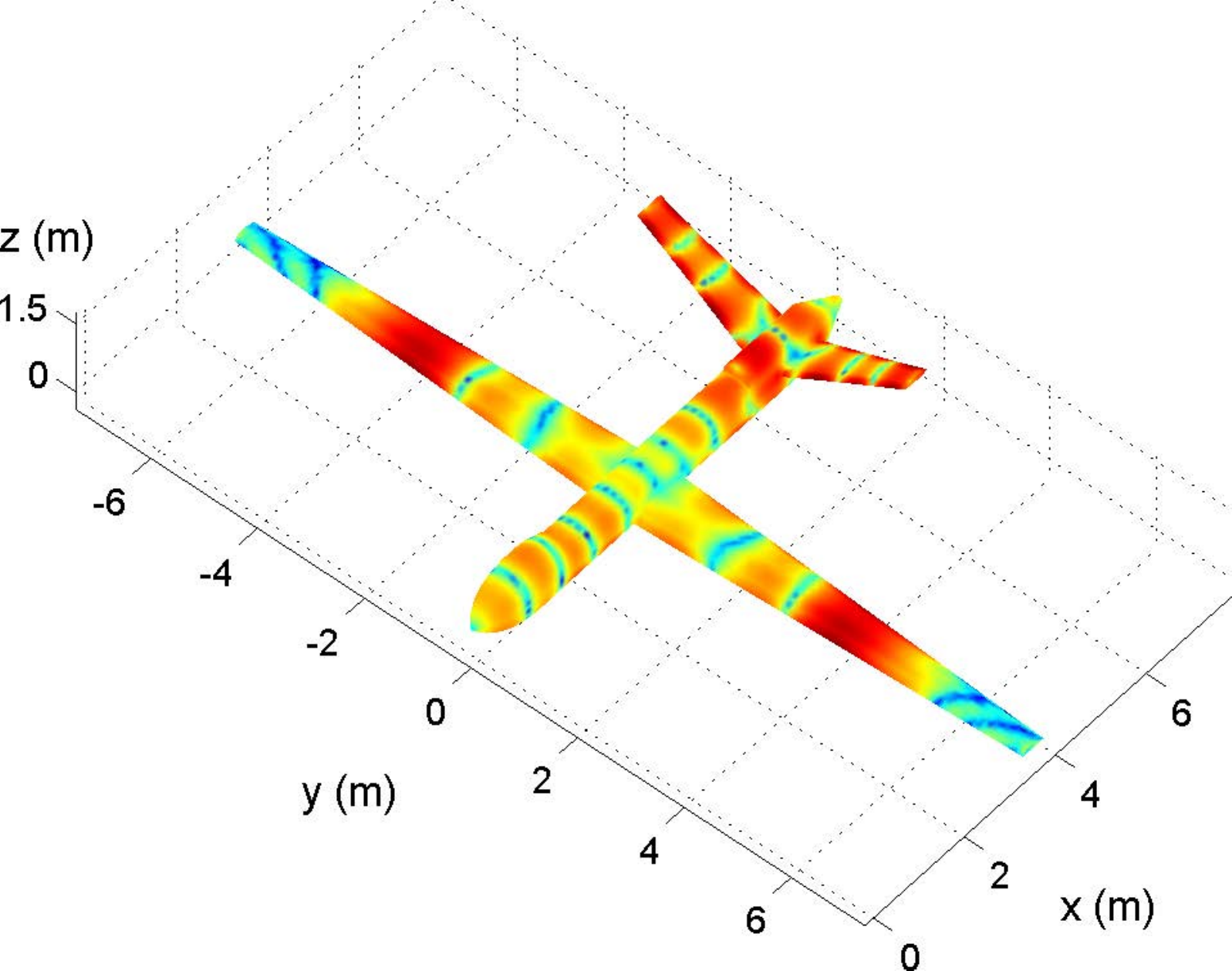} \label{CM4Pred}} 
\caption{Characteristic modes of UAV computed using MLFMA + IRA ($21\thinspace056$ unknowns): (a) Mode with characteristic value $\lambda = -0.0548$. (b) Mode with characteristic value $\lambda = 0.0123$.}  \label{CM Pred}
\end{figure}

We finally consider large platforms, i.e., aerospace vehicles, where we can only apply the MLFMA + IRA scheme for CM analysis. JDQZ becomes inefficient in such examples as there are lack of preconditioners for expediting the convergence to solutions with reasonable accuracy when (\ref{JDCE}) is solved. We compute several CMs of an unmanned aerial vehicle (UAV), of which the length, width and height are $8.06$~m, $14.78$~m, and $2.16$~m,  or $5.72\lambda$, $10.49\lambda$, and $1.53\lambda$ at the frequency of $213$~MHz, respectively. In this simulation, the number of unknowns is $21\thinspace056$. The GMRES tolerance is set to $5\times10^{-3}$, and $\epsilon_1$, $\epsilon_2$ and $\epsilon_3$ in SAI are set to $0.01$, $0.012$ and $0.07$, respectively. The computed current patterns (magnitudes) are illustrated in Fig.~{\ref{CM Pred}}.

\begin{figure}[!t]
\centering
\subfigure[]{\includegraphics[scale=0.33]{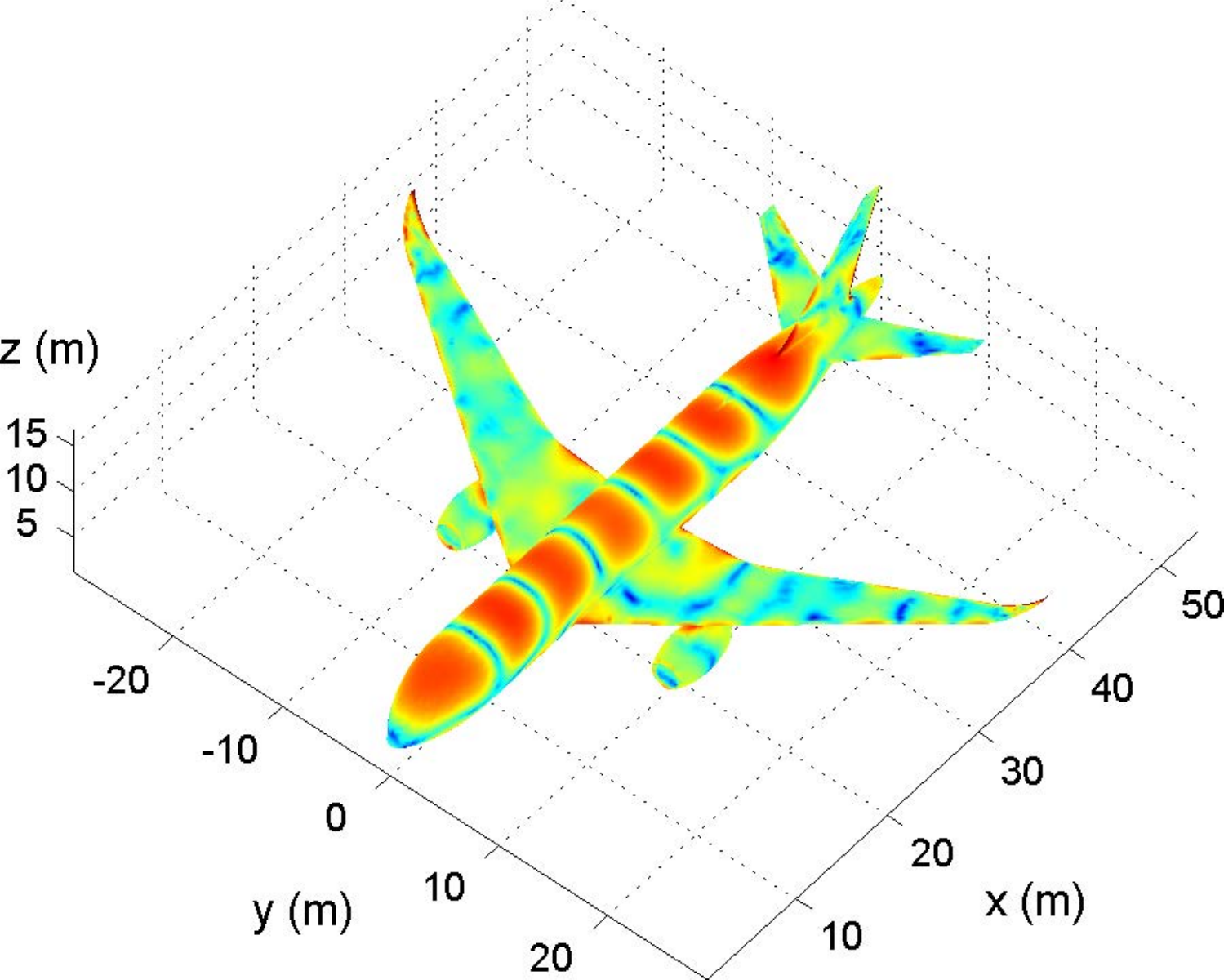} \label{CM1Pred}}
\subfigure[]{\includegraphics[scale=0.33]{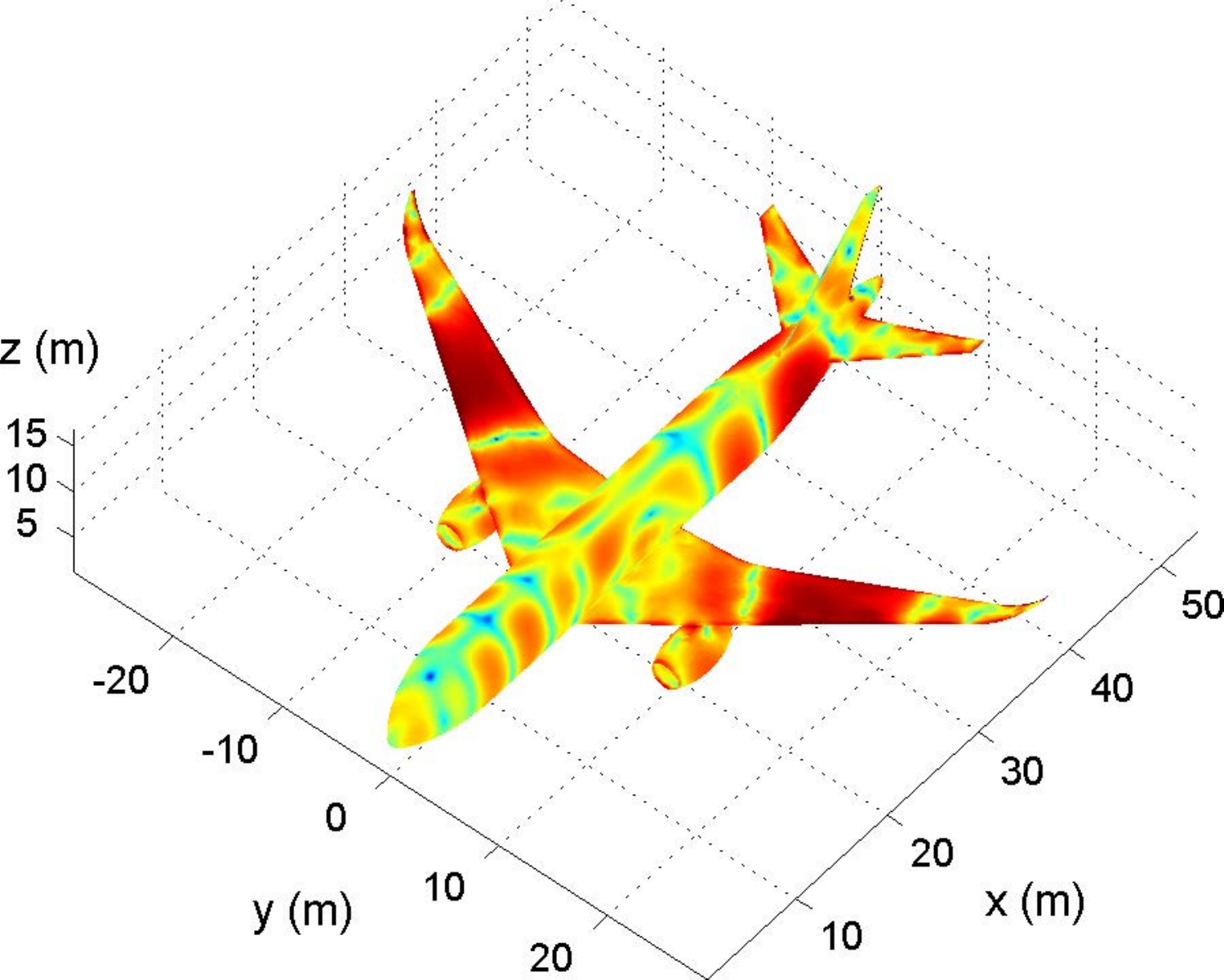} \label{CM4Pred}} 
\caption{Characteristic modes of Boeing 787 Dreamliner computed using MLFMA + IRA ($117\thinspace834$ unknowns): (a) Mode with characteristic value $\lambda =0.0358$. (b) Mode with characteristic value $\lambda = -0.0424$.}  \label{CM B787}
\end{figure}

The CM analysis is next performed on a Boeing 787 Dreamliner. The length, width, and height of the plane is $53.83$~m, $58.41$~m, and $15.70$~m, respectively, which are $7.72\lambda$, $8.37\lambda$, and $2.25\lambda$ at $43$~MHz. Totally $117\thinspace834$ unknowns are used in this simulation. The GMRES tolerance is again set to $5\times 10^{-3}$, and $\epsilon_1$, $\epsilon_2$ and $\epsilon_3$ in SAI are set to $0.008$, $0.01$ and $0.07$, respectively. Certain characteristic currents of the Dreamliner are depicted in Fig~{\ref{CM B787}}. The color axis is adjusted to clearly manifest the current patterns on the plane body. 

\begin{table}
\begin{center}
 \caption{Computational Detail in Different Numerical Examples Using MLFMA + IRA}
 \label{TCD}    
\begin{tabular}{ccccccc} 
 \hline
 Example & nev  &  ncv & nlv & nii & Memory & CPU time \\
 \hline
 Plate    & 5  &  20    &  4  &  $>$45 & 20~MB & 15 secs \\
 Sphere & 10 &  120  &  5  & $>$170  & 300~MB & 55 mins \\
 UAV    & 4   &  140  &   7  & $>100$ & 436~MB & 2.1 hrs\\
 Dreamliner  & 4    & 180  &  7  &  $>240$  & 3.1~GB & 18 hrs\\
 \hline
 \end{tabular}
 \end{center}
 $^*$ Note that in the examples of UAV and Dreamliner, one may obtain $4$ converged eigenpairs with smaller ncv. In the sphere example, $15\thinspace015$ unknowns are used.
\end{table}

Numerical examples in this study are computed with a fundamental, non-parallel MLFMA program and the open-source ARPACK on an Intel Core i5-2400 CPU with $3.10$ GHz clock rate. Computational detail in different examples are summarized in Table~\ref{TCD}, where nlv is the total number of MLFMA levels (including level-$0$), nev is the number of eigenvalues needed, ncv is the maximum number of Arnoldi vectors used, and nii is the number of inner iterations for calculating each $\bar{\mathbf{Z}}^{-1} \mathbf{u}$. It can be inferred that the SAI preconditioner becomes less efficient as the problem scale grows large. For highly ill-conditioned impedance matrix, it calls for more efficient preconditioning schemes such as the Calderon projector to further enhance the performance of the proposed CM analysis. 

\section{Conclusion}
We develop an MLFMA based CM analysis for large scale applications which severely challenge the conventional MoM based approach. The standard MLFMA program can be easily modified  to perform the required matrix-vector product operations in the iterative eigensolvers.  Numerical examples are provided to demonstrate the validity and efficiency of the proposed scheme. In even larger scale simulations, MLFMA can be implemented in parallel on many processors for better performance. The same idea can be easily extended to an incorporation of the mixed-form fast multipole algorithm (MF-FMA) \cite{Jiang05} and CM analysis, which is beneficial for models with fine structures.

\section*{Acknowledgment}
%
%
The authors would like to thank Q. S. Liu for her helpful discussion.

\ifCLASSOPTIONcaptionsoff
  \newpage
\fi

\bibliographystyle{IEEEtran}
\bibliography{thesisrefs}







\end{document}